\documentclass[nosumlimits]{amsart}
\usepackage{amscd,amssymb}

\theoremstyle{plain}
\newtheorem{prop}[subsection]{Proposition}
\newtheorem{thm}[subsection]{Theorem}
\newtheorem{lem}[subsection]{Lemma}
\newtheorem{cor}[subsection]{Corollary}
\newtheorem*{theo}{Theorem}

\theoremstyle{definition}
\newtheorem{exm}[subsection]{Example}
\newtheorem{rem}[subsection]{Remark}

\numberwithin{equation}{section}

{%
 \begin{enumerate}
}
{%
 \end{enumerate}%
}

\newenvironment{romenum}%
{%
 \begin{enumerate}
}
{%
 \end{enumerate}%
}

\newcommand{\A}{{\mathcal A}}
\newcommand{\B}{{\mathcal B}}
\newcommand{\CC}{{\mathcal C}}
\newcommand{\RR}{{\mathcal R}}

\newcommand{\W}{{\mathcal W}}

\newcommand{\cV}{{\mathcal V}}

\newcommand{\Z}{{\mathbb Z}}

\newcommand{\C}{{\mathbb C}}
\newcommand{\CP}{{\mathbb{CP}}}
\newcommand{\T}{{({\mathbb C}^*)^n}}
\newcommand{\F}{{\mathbb F}}

\newcommand{\bT}{{\mathbf T}}
\newcommand{\V}{{\mathbf V}}
\newcommand{\bone}{{\mathbf 1}}
\newcommand{\bz}{{\mathbf 0}}
\newcommand{\bt}{{\mathbf t}}

\newcommand{\dA}{{\mathbf{d}\mathcal{A}}}

\newcommand{\cA}{{\mathbf{c}\mathcal{A}}}

\newcommand{\bl}{{\boldsymbol{\lambda}}}
\newcommand{\bD}{{\boldsymbol{\Delta}}}
\newcommand{\bp}{{\boldsymbol{\partial}}}
\newcommand{\bd}{{\boldsymbol{\delta}}}

\newcommand{\p}{\partial }
\newcommand{\M}{{M'}}

\newcommand{\ii}{\text{i}}

\newcommand{\locB}{{B^{\operatorname{cc}}}}
\newcommand{\loc}{{\operatorname{cc}}}

\renewcommand{\a}{{\alpha }}
\renewcommand{\b}{{\beta }}
\renewcommand{\c}{{\gamma }}
\renewcommand{\d}{{\delta }}
\renewcommand{\l}{{\lambda}}
\renewcommand{\L}{{\Lambda }}
\renewcommand{\ll}{{\ell }}
\renewcommand{\*}{{\bullet}}

\DeclareMathOperator{\rank}{rank}

\DeclareMathOperator{\coker}{coker}
\DeclareMathOperator{\id}{id}
\DeclareMathOperator{\Aut}{Aut}

\DeclareMathOperator{\Mat}{Mat}

\DeclareMathOperator{\GL}{GL}
\DeclareMathOperator{\ab}{ab}

\DeclareMathOperator{\pr}{pr}

\DeclareMathOperator{\TT}{T}

\DeclareMathOperator{\Bb}{B}
\DeclareMathOperator{\Dd}{D}

\begin{document}

\title[Characteristic Varieties of Arrangements]
{Characteristic Varieties of Arrangements}

\author[Daniel C.~Cohen]
{Daniel C.~Cohen$^1$}
\address{Department of Mathematics,
Louisiana State University,
Baton Rouge, LA 70803}
\email{cohen@math.lsu.edu}
\urladdr{http://math.lsu.edu/\~{}cohen}

\author[Alexander I.~Suciu]
{Alexander I.~Suciu$^2$}
\address{Department of Mathematics,
Northeastern University,
Boston, MA 02115}
\email{alexsuciu@neu.edu}
\urladdr{http://www.math.neu.edu/\~{}suciu}


\thanks{$^1$ Partially supported by
grant LEQSF(1996-99)-RD-A-04 from the Louisiana Board of Regents.}
\thanks{$^2$ Partially supported by NSF~grant
DMS--9504833.}



\begin{abstract} 
The $k^{\text{th}}$ Fitting ideal of the Alexander invariant $B$
of an arrangement $\A$ of $n$ complex hyperplanes defines a
characteristic subvariety, $V_k(\A)$, of the algebraic 
torus $\T$.  In the combinatorially determined case where $B$ 
decomposes as a direct sum of local Alexander invariants, we 
obtain a complete description of $V_k(\A)$.  For any arrangement 
$\A$, we show that the tangent cone at the identity of this variety  
coincides with $\RR^{1}_{k}(A)$, one of the cohomology support 
loci of the Orlik-Solomon algebra.  Using work of Arapura~\cite{Ar} 
and Libgober~\cite{L3}, we conclude that all positive-dimensional 
components of $V_{k}(\A)$ are combinatorially determined, and that 
$\RR^{1}_{k}(A)$ is the union of a subspace arrangement 
in $\C^n$, thereby resolving a conjecture of Falk~\cite{Fa}.  
We use these results to study the reflection arrangements 
associated to monomial groups.
\end{abstract}

\maketitle

\section*{Introduction}
\label{sec:intro}

A hyperplane arrangement is a finite collection $\A$ of 
codimension one subspaces in a finite-dimensional complex 
vector space $V$. Two principal objects associated to $\A$ are   
the complement, $M(\A)=V\setminus \bigcup_{H\in \A}\, H$, 
and the intersection lattice,  
$L(\A)=\{\bigcap_{H\in \B}\, H \mid \B\subseteq \A\}$.  
A central problem in the study of arrangements 
is to elucidate the relationship 
between these two seemingly disparate objects---one topological, 
the other combinatorial.  The paradigmatic result in 
this direction is the theorem of Orlik and Solomon~\cite{OS}, 
which asserts that the cohomology ring of the complement, 
$H^{*}(M(\A);\C)$, is isomorphic to a certain algebra, $A(\A)$, 
which is completely determined by the lattice $L(\A)$.  

The above result leads one to investigate the extent to which
the lattice $L=L(\A)$ determines the topology of the complement $M=M(\A)$.  
Examples of Rybnikov~\cite{Ryb} show that the fundamental group, 
$G=\pi_1(M)$, is not combinatorially determined in general.  
However, certain invariants of $G$, such as the lower central 
series quotients (cf. Falk~\cite{falkLCS}), are determined by 
the lattice.  Thus, it is natural to ask whether a given 
isomorphism type invariant of $G$ is determined by the 
isomorphism type of $L$.  In this paper, we show that the  
central characteristic subvarieties of $G$ are indeed 
combinatorially determined.

Let $M'$ be the maximal abelian cover of $M$, and let 
$B(\A)=H_{1}(\M;\Z)$ be the Alexander invariant of $\A$, viewed 
as a module over the Laurent polynomial ring $\L=\Z\Z^n$, where 
$n=|\A|$.  From the presentation of $B=B(\A)$ found in \cite{CSai}, 
several invariants of the group $G=\pi_1(M)$ may be computed.  
One such invariant is the $k^{\text{th}}$~Fitting ideal, 
$F_{k}(B)$.  This is the ideal of $\L$ generated by all the 
codimension~$k-1$ minors of a presentation matrix for $B$, 
and is well-known to be independent of the presentation.  
Let $\T$ be the complex algebraic $n$-torus, with coordinate 
ring $\L_{\C}=\C\Z^n$.  The $k^{\text{th}}$ 
characteristic variety of $\A$ is the subvariety $V_k(\A)$ of $\T$ 
defined by the ideal $F_k(B)\otimes \C$ of $\L_{\C}$.   
These varieties depend only on the group $G$, 
up to a monomial change of basis in $\T$.  
The characteristic varieties first appeared in a more general context in 
\cite{Nov}, \cite{DF}, \cite{L1}, and have been recently studied 
in \cite{Ar}, \cite{Hir}, \cite{L3}.  

The lattice $L$ of $\A$ is a partially ordered set, ordered by 
reverse inclusion, with rank function given by codimension.  
(See Orlik and Terao~\cite{OT} as a general reference for arrangements.)
Each flat $X\in L_2$ of rank two gives rise to a subvariety $V_X$ of 
$V_k(\A)$.  In the combinatorially determined instance (identified 
in \cite{CSai}) where the module $B(\A)$ decomposes as a direct sum 
of ``local'' Alexander invariants, these subvarieties $V_X$ provide 
a complete description of the characteristic varieties of $\A$.

In general, the characteristic varieties of an arrangement $\A$ possess 
``non-local'' irreducible components.  To analyze such components, we make 
use of the relation between these varieties and the cohomology of 
local systems on the complement $M$ of $\A$.  Each point $\bt\in\T$ 
gives rise to a local coefficient system $\C_\bt$ on $M$.  The 
characteristic variety $V_{k}(\A)$ may be identified 
with the cohomology support locus 
$W^1_{k}(M)=\{ \bt \in \T \mid \rank H^1(M;\C_\bt) \ge k\}$.  This 
variety was studied by Arapura~\cite{Ar}, who showed that it is a 
union of torsion-translated subtori in $\T$ in more general circumstances.  
In his work on characteristic varieties of plane algebraic 
curves~\cite{L2}, \cite{L3}, 
Libgober proved that all positive-dimensional components of $V_k(\A)$ 
pass through the identity element $\bone=(1,\dots ,1)$ of $\T$.

A comparable analysis in the Orlik-Solomon algebra has recently been 
carried out by Falk~\cite{Fa}.  Let $A=A(\A)$ be the OS-algebra of $\A$, 
generated by $a_1,\dots,a_n$.  Each point $\bl \in \C^n$ determines an 
element $\omega=\sum_{i=1}^n \l_i a_i$ of $A^1$.  Falk's invariant 
varieties are the cohomology support loci 
$\RR^{1}_{k}(A)=\{\bl \in \C^n \mid \rank H^1(A,\mu)\ge k\}$, 
where $A$ is viewed as a complex with differential 
$\mu(\eta)=\omega \wedge \eta$.  We establish the relation between 
these varieties and the characteristic varieties $V_{k}(\A)$.  We prove:
\begin{theo} 
\label{th:tconecv}
The tangent cone $\cV_{k}(\A)$ of $V_{k}(\A)$ at the point $\bone$ 
coincides with $\RR^{1}_{k}(A)$.
\end{theo}
This result has several significant consequences.
First, together with Arapura's work noted above, it shows
that $\RR^{1}_{k}(A)$ is the union of a subspace arrangement in $\C^n$, 
resolving a conjecture of Falk~\cite{Fa}.
Second, in conjunction with Libgober's result, it allows
us to conclude that the monomial isomorphism type of the 
central characteristic subvariety $\check{V}_k(\A)$ 
(the subvariety of $V_{k}(\A)$ consisting of 
those irreducible components passing through $\bone$) 
is determined by the isomorphism type of the lattice $L(\A)$.
Finally, the identification $\RR^{1}_{k}(A) =\cV_{k}(\A)$ provides
a combinatorial means for detecting non-local components of $V_{k}(\A)$.  
We illustrate this technique by analyzing the first characteristic
varieties of braid arrangements and monomial arrangements.
As a consequence, we show that the associated generalized
pure braid groups are not isomorphic to the corresponding
products of free groups, extending a result of~\cite{CSpn}.

\section{Fitting Ideals and Characteristic Varieties}
\label{sec:DetIdCharVar}

In this section we review the definition of the Alexander
invariant, and the associated Fitting ideals and characteristic 
varieties, for a finite complex.

\subsection{Alexander Invariant}
\label{subsec:Bmodule}
Let $M$ be a path-connected space that has the homotopy type of a 
finite CW-complex.  Let $G=\pi _{1}(M,*)$ be the fundamental group, 
and $K=H_{1}(M)$ its abelianization.  Let $\M$ be the maximal 
abelian cover of $M$, with group of deck-transformations identified 
with $K$.  The action of $K$ on $M$ induces an action 
on the homology groups $H_*(M')$. This defines on $H_*(M')$ the 
structure of a module over the group ring $\Z K$.  The $\Z K$-module 
$B=H_{1}(M')$ is called the (first) {\em Alexander invariant} of $M$.  
This module depends only on the group $G$.  Indeed, $B=G'/G''$, 
with the action of $K=G/G'$ defined by the extension 
$1\to G'/G''\to G/G'' \to G/G'\to 1$.  

Since $M$ is up to homotopy a finite complex, $G$ is a finitely
presented group.  Now assume $K$ is free abelian, and choose a system of
generators, $t_{1},\dots , t_{n}$.   This provides an identification 
of the group ring $\Z K$ with the ring of finite  Laurent series 
in $n$ (commuting) variables, $\L =\Z [t_{1}^{\pm 1},\dots ,t_{n}^{\pm 1}]$.  
Thus, we may view the Alexander invariant $B$ as a module over $\L$.  
As was shown by Crowell, $B$ is in fact a finitely presented 
$\L$-module.  See \cite{Hil} for details.

\subsection{Fitting Ideals} 
\label{sec:detideals}  
Recall a standard notion from commutative algebra.  
Let $A$ be a module over a commutative ring $R$, 
and assume $A$ has free presentation 
$R^p\xrightarrow{\Omega} R^q\to A\to 0$.  
The {\em $k^{\text{th}}$ Fitting ideal} $F_k(A)$ of 
the module $A$ is the ideal of $R$ generated by 
the $(q-k+1)\times (q-k+1)$ minors of the matrix $\Omega$.  
We set $F_k(A)=0$ if $k\le 0$ or $k\le q-p$, and $F_k(A)=R$ 
if $k>q$.  The ideal $F_k(A)$---also known 
as the $(k-1)^{\text{th}}$ determinantal (or elementary) 
ideal of $A$---is independent of the choice of presentation 
for $A$, see e.g.~\cite{No}. The Fitting ideals form an 
ascending chain $0=F_{0}\subseteq F_{1}\subseteq 
\cdots \subseteq F_{q}\subseteq F_{q+1}=R$.   

Now let $M$ be a finite CW-complex, with $H_1(M;\Z)=\Z^n$.   
The Fitting ideals, $F_k(M)$, of the Alexander invariant 
$H_1(M')$ are homotopy-type invariants of $M$.  
More precisely, if $f:M \to N$ is a homotopy equivalence, 
the extension of $f_*:H_1(M;\Z) \to H_1(N;\Z)$ to group rings 
restricts to an isomorphism $f_*:  F_k(M)\to F_k(N)$ between 
the corresponding ideals, see~\cite {Fox}, \cite{Tu}.  

\subsection{Characteristic Varieties} 
\label{sec:charvars} 
Let $\L_{\C}=\L\otimes \C$ be the ring of Laurent polynomials with 
complex coefficients. This is the coordinate ring of the algebraic 
torus $\T$.  Consider the complexified Alexander invariant of $M$, 
$B_{\C}(M)=H_1(M'; \C)$, viewed as a module over $\L_{\C}$, and  
let $F_{k}^{\C}(M)$ be its $k^{\text{th}}$ Fitting ideal.
Following Libgober~\cite{L1}, we call 
the reduced algebraic variety $V_k(M)=\V(F_{k}^{\C}(M))$ 
the {\em $k^{\text{th}}$ characteristic variety} of $M$. 
The characteristic varieties form a descending tower 
$\T=V_{0}\supseteq V_{1}\supseteq \cdots \supseteq V_{q}
\supseteq V_{q+1}=\emptyset$.

The varieties $V_k(M)$ are homotopy-type invariants of $M$.   
More precisely, if $M$ is homotopy equivalent to $N$, 
there exists a monomial automorphism $g:\T\to \T$, given by 
$g(t_i)=t_1^{a_{i,1}} \cdots t_n^{a_{i,n}}$, for some matrix 
$(a_{i,j})\in\GL(n,\Z)$, such that $g(V_k(M))=V_k(N)$.  

The Fitting ideals and characteristic varieties depend 
only on the fundamental group $G=\pi_1(M)$, up to a 
change of basis as above. Accordingly, we shall denote 
them, when convenient, by $F_k(G)$, resp.~$V_k(G)$.  
When $G$ is the group of an arrangement $\A$, we may
use the notation $F_k(\A)$, resp.~$V_k(\A)$.

\begin{exm}
\label{exm:Fn} 
Let $G = \F_n $ be the free group of rank $n$.  Recall the 
standard free resolution $C_{\bullet}$ of $\Z$ over $\L=\Z\Z^n$:
\begin{equation}  
\label{eq:znres}
0\rightarrow C_{n} \xrightarrow{d_{n}}
\cdots \rightarrow C_{3} \xrightarrow{d_{3}}
C_{2} \xrightarrow{d_{2}}
C_{1} \xrightarrow{d_{1}}
C_{0} \xrightarrow{\epsilon}
\Z  \rightarrow  0, 
\end{equation}
where $C_{0}=\L $, $C_{1}=\L ^{n}$, 
and $C_{k}=\bigwedge ^{k} C_{1}=\L^{\binom {n}k}$, 
and the differentials are given by 
$d_{k}(e_J)=\sum _{r=1}^{k} (-1)^{r} (t_{j_{r}}-1)
\cdot e_{J\setminus\{j_r\}}$, where
$e_J=e_{j_1}\wedge\dots\wedge e_{j_k}$ if
$J=\{j_1,\dots,j_k\}$. 

Let $B=B(\F_n)$ be the Alexander invariant, 
$F_k=F_k(\F_n)\subset \Lambda$ the Fitting ideals, and 
$V_k=V_k(\F_n)\subset \T$ the characteristic varieties.  
If $n=1$, we have $B=0$, and so $F_0=0$, $F_1=\Lambda$, and 
$V_0=\C^*$, $V_1=\emptyset$.  If $n\ge 2$, a free presentation 
for $B$ is given by $C_3 \xrightarrow{d_3} C_2\to B \to 0$.  
A standard computation yields:
\[
F_k=
\begin{cases}
0&\text{for }\, 0\le k \le n-1,\\
I^{\binom{n}{2}-k+1} &\text{for }\, n\le k\le\binom{n}{2},\\
\Lambda &\text{for }\, k>\binom{n}{2},
\end{cases}
\qquad
V_k=
\begin{cases}
\T&\text{for }\, 0\le k \le n-1,\\
\bone&\text{for }\, n\le k\le\binom{n}{2},\\
\emptyset &\text{for }\, k>\binom{n}{2},
\end{cases}
\]
where $I=(t_1-1,\dots ,t_n-1)$ is the augmentation ideal of $\L$ 
and $\bone=(1,\dots ,1)$ is the identity element of $\T$. 
\end{exm}

\section{Braid Monodromy and Alexander Invariant of an Arrangement}
\label{sec:BmonoAinv}

In this section, we review the algorithms for determining the  
braid monodromy and the Alexander invariant of the complement 
of a complex hyperplane arrangement, as developed  
in \cite{CSbm}, \cite{CSai}, and record some immediate consequences.

\subsection{Braid Monodromy} 
\label{subsec:bmono}
The fundamental group of the complement of a complex hyperplane 
arrangement is isomorphic to that of a generic two-dimensional 
section.  So, for the purpose of studying the Alexander invariant, 
it is enough to consider affine line arrangements in $\C ^{2}$.
Let $\A =\{H_{1},\dots ,H_{n}\}$ be such an arrangement,
with complement $M=M(\A)=\C^2\setminus \bigcup_{i=1}^{n} H_i$, and 
vertices $\mathcal V = \{v_{1},\dots ,v_{s}\}$.
If $v_k = H_{i_1}\cap \dots \cap H_{i_r}$, let $X_k=\{i_1,\dots,i_r\}$
denote the corresponding ``vertex set.''  We identify the
set $L_2=L_2(\A)$ of rank two elements in the lattice of $\A$ and the
collection $\{X_1,\dots, X_s\}$ of vertex sets of $\A$.

The braid monodromy is determined as follows (see~\cite{CSbm} 
for details).  Choose coordinates $(x,z)$ in $\C ^{2}$ so that the 
first-coordinate projection map is generic with respect to $\A $.   
Let $f(x,z)=\prod _{i=1}^{n} (z-a_{i}(x))$ be a defining polynomial 
for $\A $.  The root map $a=(a_{1},\dots , a_{n}):\C \to \C ^{n}$ 
restricts to a map from the complement of $\mathcal{Y}=
\pr _{1}(\mathcal{V})$ to the complement of the braid arrangement 
$\A _n=\{\ker (y_{i}-y_{j})\}_{1\le i <j \le n}$.  Identify 
$\pi_{1}(\C \setminus \mathcal{Y})$ with the free group $\F_{s}$, 
and the group of the braid arrangement, $\pi_1(M(\A_n))$, 
with the pure braid group $P_{n}$.  Then, the {\em braid monodromy} 
of $\A$ is the homomorphism on fundamental groups,
$\alpha:\F_s \to P_n$, induced by the root map.

The generators, $\{\a_{1},\dots, \a_s\}$, of the image of 
the braid monodromy can be written explicitly 
using a {\em braided wiring diagram} $\W$ associated to $\A$.
Such a diagram, determined by the choices made above, may be 
specified by a sequence of vertex sets and braids,
$\W =\W_s= \{ X_1,\b_1,X_2,\b_2,\dots$, $\b_{s-1},X_s\}$.
The braid monodromy generators are given by
$\a _{k}=A_{X_{k}}^{\d _{k}}$, where
$A_{X_k}$ is the full twist on the strands comprising $X_k$, and
$\d_k$ is a pure braid determined by the subdiagram $\W_k$.

A presentation for the fundamental group of 
$M$ may be obtained from the braid monodromy generators 
using the Artin representation $P_n\to \Aut(\F_n)$:
\begin{equation*}
\label{eq:pi1}
\pi_1(M)=\langle \c_1,\dots ,\c_n \mid \a_k(\c_i) = \c_i \quad 
\text{for } i=1,\dots ,n\, \text{ and }\, k=1,\dots ,s\rangle .
\end{equation*}
Note that all relations are commutators.  This presentation can be
simplified by Tietze-II moves---eliminating redundant relations.
For each $X\in L_2$, let $X'=X\setminus\{\min X\}$.  Then, the
braid monodromy presentation is
\begin{equation}  
\label{eq:bmpres}
\pi_1(M)=\langle \c_1,\dots ,\c_n \mid \a_k(\c_i) = \c_i \quad  
\text{for } i \in X_k'\, \text{ and }\, k=1,\dots ,s\rangle .
\end{equation}
The number of relations in this presentation is 
equal to the second Betti number,  
$b= b_2(M)=\sum_{X\in L_2} |X'|$, 
of the complement.

\subsection{Alexander Invariant} 
\label{subsec:AlexInv}
For an endomorphism $\a$ of the free group $\F_{n}$, let
$\Theta(\a):C_{1}\to C_{1}$ be its abelianized Fox Jacobian.  
This is a $\L$-linear map, whose matrix has rows 
\[
\Theta (\a )(e_{i})=\nabla ^{\ab }(\alpha (t_{i}))=
\sum _{i=1}^{n} \left( \frac{\p\alpha (t_{i})}{\p t_{i}}\right)^{\ab}
e_{i}.
\]
The restriction of $\Theta$ to the pure braid group $P_n<\Aut(\F_n)$ 
is the Gassner representation, $\Theta:P_n\to \GL(n,\L)$.  
Applying the Fox Calculus to the presentation \eqref{eq:bmpres}, 
we obtain the chain complex, $C_\*(M')$, of the maximal abelian 
cover of $M$:
\begin{equation} 
\label{eq:m'cx}
\L^b \xrightarrow{\p_{2}}
\L^n \xrightarrow{\p_{1}}
\L \xrightarrow{\epsilon}
\Z  \rightarrow  0, 
\end{equation}
where $\p_1=d_1=\pmatrix t_1-1&\cdots&t_n-1 \endpmatrix^\top$, 
and $\p_2$ is the Alexander matrix, with rows 
$(\Theta(\a_k)-\id)(e_i)$, indexed by 
$i \in X_k'$ and $k\in\{1,\dots ,s\}$.  
If $X$ is a non-empty subset of $[n]=\{1,\dots,n\}$, 
let $C_1[X]$ denote the submodule of $C_1=\L^n$
spanned by  $\{ e_j \mid j \in X\}$. 
Then, via the identification $C_2(M') \cong \bigoplus_{X \in L_2} C_1[X']$, 
the Alexander matrix may be viewed as a map
$\p_2: \bigoplus_{X \in L_2} C_1[X'] \to C_1$.

A presentation for the Alexander invariant, $B=H_1(M')$, may be
obtained by comparing the complex \eqref{eq:m'cx} with the 
resolution \eqref{eq:znres}.  We paraphrase the construction 
of~\cite{CSai}.  For $\a\in P_n$, with Artin representation 
given by $\a(t_{i})=z_{i}t_{i}z_{i}^{-1}$, the homomorphism 
$\Phi(\a):C_{1}\to C_{2}$ defined by
\begin{equation*}  \label{eq:Phidef}
\Phi(\a)(e_{i})=e_{i}\wedge \nabla ^{\ab}(z_{i})
\end{equation*}
satisfies $\Theta(\a)-\id=d_{2}\circ \Phi (\a)$.  
If $\a=A_X$ is the full twist on $X=\{i_{1},\dots ,i_{r}\}$, 
given in terms of the standard generators of $P_n$ by 
\[
A_{X} =
(A_{i_{1},i_{2}})(A_{i_{1},i_{3}}A_{i_{2},i_{3}})(A_{i_{1},i_{4}}
A_{i_{2},i_{4}}A_{i_{3},i_{4}})\cdots
(A_{i_{1},i_{r}}\cdots A_{i_{r-1},i_{r}}), 
\]
then $\Phi (A_{X}):C_{1} \to C_{2}$ is given by
\begin{equation*}
\label{eq:PhiX}
\Phi(A_X)(e_{i})
    =\begin{cases}
      e_{i}\wedge \nabla _{X}  &\text{if } i\in X\\
      (t_{i}-1)\nabla _{X}\wedge \nabla _{^{i}X}
      &\text{if } i_{1}\le i \le i_{r} \text{ and } i\notin X \\
      0 &\text{otherwise}, 
      \end{cases}
\end{equation*}
where 
$X^{i}=\{ j\in X \mid j<i\}$, $^{i}X=\{j\in X\mid i<j\}$, 
$t_{X}^{}=\prod _{j\in X} t_{j}$, and 
$\nabla_X=\sum_{j\in X}t_{X^{j}}e_{j}$.  
In general, if $\a=A_X^\d$, then 
$\Phi(\a)=\Theta_{2}(\d)\circ\Phi(A_X)\circ\Theta (\d^{-1})$, 
where $\Theta_{k}=\bigwedge^{k}\Theta:C_{k}\to C_{k}$.

If $A_X^\d$ is the braid monodromy generator corresponding to 
$X\in L_2$, let $\Phi_X$ denote the restriction of $\Phi(A_X^\d)$ 
to $C_1[X']$.  Then the map $\Phi:\bigoplus_{X \in L_2} C_1[X'] \to C_2$, 
defined by $\Phi|_{C_1[X']}=\Phi_X$, satisfies
\begin{equation} 
\label{eq:phid2}
\p_2 = d_2 \circ \Phi.
\end{equation}
Thus, the maps $\id_{C_{0}}$, $\id_{C_{1}}$, and $\Phi$ constitute 
a chain map $\Phi_\*$ from the complex \eqref{eq:m'cx} to the resolution 
\eqref{eq:znres}.  

\begin{thm} \label{thm:aipresthm}
The Alexander invariant of the arrangement $\A$ has presentation
\begin{equation*} 
\label{eq:aipres}
K_1 \xrightarrow{\Delta} K_0 \to B \to 0,
\end{equation*}
where $K_0=C_2$, 
$K_1=C_2(M') \oplus C_3 = \bigoplus_{X \in L_2} C_1[X'] \oplus C_3$, 
and $\Delta = \bigl(\smallmatrix \Phi \\ d_3 \endsmallmatrix\bigr)$.
\end{thm}
\noindent This presentation has $\binom{n}{2}$ generators and
$b+\binom{n}{3}$ relations, where $b=b_2(M)$.
\begin{proof}
Let $K_\*(\Phi)$ denote the mapping cone of the chain map $\Phi_\*$.  
We then have an exact sequence of complexes, 
$0 \to C_\* \to K_\*(\Phi) \to C_{\*-1}(M') \to 0$, explicitly:
\begin{equation}
\label{eq:mapconediagram}
\begin{CD}
@>>>  C_3 @>d_3>> C_2 @>d_2>> C_1 @>d_1>> C_0 \\
&& @VVV  @VVV 
@VVV  @VVV \\
@>>>  {C}_{2}(M') \oplus C_3   
@>{\bigl(\smallmatrix \p_2 & \Phi \\ 0 & d_3 \endsmallmatrix\bigr)}>> 
C_1 \oplus C_2
@>{\bigl(\smallmatrix d_1 & -\id \\ 0 & \ d_2 \endsmallmatrix\bigr)}>> 
C_0 \oplus C_1 
@>{\bigl(\smallmatrix \id \\ d_1 \endsmallmatrix\bigr)}>> C_0 \\
&& @VVV  @VVV  @VVV \\
&& {C}_{2}(M') @>\p_2>> C_1 @>d_1>> C_0 \\
\end{CD}
\end{equation}
where the chain maps are the natural inclusion and projection.  Since 
$(C_\*,d_\*)$ is a resolution, we have $B=H_1(M') = H_2(K_\*(\Phi))$, 
and it is immediate from the above diagram that $H_2(K_\*(\Phi))=\coker
\Delta$.
\end{proof}

\subsection{Characteristic Varieties}
\label{subsec:cva}
We now describe the characteristic varieties of the arrangement $\A$ in 
terms of the presentation of the Alexander invariant obtained above.  
First we establish some notation.

Consider the evaluation map $\L\times \T\to \C$, 
which takes a Laurent polynomial in $n$ variables, $f$, and a point 
$\bt=(t_1,\dots ,t_n)$ and yields $f(\bt)=f(t_1,\dots ,t_n)$.  
For fixed $f\in\L$, we get a map $\mathbf{f}:\T\to\C$.  
More generally, we have the map 
$\Mat_{p\times q}(\L)\times \T\to \Mat_{p\times q}(\C)$, which 
takes a matrix $F: \L^p\to \L^q$ and evaluates each entry at $\bt$, 
to give $F(\bt): \C^p\to \C^q$. For fixed $F\in\Mat_{p\times q}(\L)$, 
we get a map $\mathbf{F}:\T\to\Mat_{p\times q}(\C)$. 

The above considerations, when applied to the presentation matrix 
$\Delta$ for the Alexander invariant from Theorem~\ref{thm:aipresthm}, yield
the following description of the characteristic varieties of $\A$:
\begin{equation}
\label{eqn:cvd}
V_k(\A) = \left\lbrace \bt \in \T \mid 
\rank \bD(\bt) \le \tbinom{n}{2}-k 
\right\rbrace.
\end{equation}

For each point $\bt \neq \bone$ in $\T$, the complex $(C_{\*}\otimes\C,
\mathbf{d}_\*(\bt))$ is acyclic.  Consequently, for each such $\bt$, we have 
$\rank \mathbf{d}_k(\bt) = \binom{n-1}{k-1}$.  In particular, 
$\rank \mathbf{d}_3(\bt)= \binom{n-1}{2}$, and it follows that 
$\rank \bD(\bt) \ge \binom{n-1}{2}$ for $\bt \neq \bone$.  For $\bt=\bone$, 
since $\mathbf{d}_3(\bone)$ is the zero matrix, we have $\rank \bD(\bone)=
\rank \boldsymbol{\Phi}(\bone)=b$ (the latter equality is a consequence of the
proof of Theorem~6.5 of \cite{CSai}, and is established by other means in the 
proof of Theorem~\ref{thm:randv} below).  Thus we have
\begin{equation}
\label{eqn:vbigk}
V_k(\A)=
\begin{cases}
\bone &\text{for $n \le k \le \binom{n}{2} - b$,}\\
\emptyset &\text{for $k > \binom{n}{2}-b$.}
\end{cases}
\end{equation}
\section{Characteristic Varieties of Decomposable Arrangements}
\label{sec:DecompVar}

In this section we study the case where the Alexander invariant 
of an arrangement decomposes as a direct sum of ``local'' invariants. 
We begin with a general formula for the characteristic varieties 
of a product of spaces.

\subsection{Products}
\label{subsec:Products}
Let $M_{1}$ and $M_{2}$ be two path-connected finite CW-complexes, 
with $K_i=H_{1}(M_{i})$ free abelian.  Let $\bT_i$ be the complex 
torus whose coordinate ring is $\C K_i$.  Let $M=M_1\times M_2$ 
be the product CW-complex, with first homology group $K=K_1\oplus 
K_2$.  
Finally, let $\bT$ be the complex torus with coordinate ring $\C K$.  

\begin{thm}\label{thm:cvprod}
With respect to the canonical decomposition $\bT=\bT_1\times \bT_2$, 
the characteristic varieties of the product $M=M_1\times M_2$ are given by
\begin{equation*} 
\label{eq:vprod}
V_k(M_1\times M_2) = (V_k(M_1)\times \bone)\cup (\bone\times 
V_k(M_2)).
\end{equation*}
\end{thm}
\begin{proof}  As noted in \cite{CSai}, the Alexander invariant 
of $M$ decomposes as $B=B_1\oplus B_2$, where 
$B_1=(H_{1}(M'_{1})\otimes _{\Z K_{1}} \Z K) \otimes _{\Z K_{2}}\Z$ 
and 
$B_2=(H_{1}(M'_{2})\otimes _{\Z K_{2}} \Z K) \otimes _{\Z K_{1}} \Z$. 
A standard property of Fitting ideals (see e.g.~\cite{No}) gives
\[
F_k(B)=\sum_{i=1}^k F_{i}(B_1)\cdot F_{k-i+1}(B_2).
\]
Thus:
\begin{align*}
V_k(M)&=\bigcap_{i=1}^k\left( \V(F_{i}(B_1))\cup
\V(F_{k-i+1}(B_2))\right)\\
&=\V(F_{k}(B_1))\cup\V(F_{k}(B_2)),
\end{align*}
and the conclusion readily follows.
\end{proof}

\begin{exm}
\label{exm:decone} 
Recall two standard constructions in arrangement theory 
(see \cite{OT} for details):  The {\em cone} of an affine 
arrangement $\A$ of $n$ hyperplanes in $\C^\ll$ is a central 
arrangement $\cA$ of $n+1$ hyperplanes in $\C^{\ll+1}$.  
The {\em decone} of a central arrangement $\A$ of $n$ hyperplanes
in $\C^\ll$ is an affine arrangement $\dA$ of $n-1$ hyperplanes
in $\C^{\ll-1}$.  The complements are related by 
$M(\A) = M(\dA) \times\C^*$.  Thus, 
$\pi_1(M(\A)) \cong \pi_1(M(\dA))\times \Z$, where the generator 
of $\Z$ can be taken to be $\c_1\cdots \c_n$, the product 
of the meridians about the hyperplanes of $\A$.  Choosing 
$H_n \in \A$ as the hyperplane ``at infinity'' in the decone 
$\dA$, we derive the following from Theorem~\ref{thm:cvprod}. 
\[
V_k(\A) = \{\bt\in \T\mid 
(t_1,\dots ,t_{n-1})\in V_k(\dA) \text{ and } t_1\cdots t_n=1\}.
\]
\end{exm}

\begin{exm} 
\label{exm:central} 
Let $\A$ be a central arrangement of $n$ lines in $\C^2$.  
Then the group of $\A$ is 
$G\cong\F_{n-1}\times \Z$, with $\Z=\langle \c_1\cdots \c_n\rangle$.  
Thus, by Example~\ref{exm:Fn} and Theorem~\ref{thm:cvprod}, 
\[
V_k(\A)=
\{\bt\in \T\mid t_1\cdots t_n=1\} \quad \text{for }\, 1\le k \le n-2,
\]
$V_k(\A)=\bone \,\text{ for }\, n-1\le k\le\binom{n-1}{2}$, and 
$V_k(\A)=\emptyset \,\text{ for }\, k>\binom{n-1}{2}$.  
\end{exm}

\begin{exm} 
\label{exm:Vfreeprod}  
More generally, the characteristic varieties of a finite direct 
product of finitely generated free groups are given by 
$V_k(\F_{n_{1}} \times \cdots \times \F_{n_{\ell}}) = 
\bigcup_{i=1}^{\ell} \widetilde V_k(\F_{n_i})$, where 
$\widetilde V_k(\F_{n_i})=
\bone\times \cdots \times V_k(\F_{n_i})\times \cdots \times \bone$
(see Example~\ref{exm:Fn}, and compare \cite{Hir}, Lemma~3.3.1). 
\end{exm}

\subsection{Local Alexander Invariants}
\label{sec:Decomp}  
Let $\A=\{H_1,\dots,H_n\}$ be an arrangement 
of $n$ lines in $\C^2$ that is transverse to infinity (that is,
no two lines of $\A$ are parallel).
For each $X\in L_2=L_2(\A)$, let $\A_X=\{H_{X,1},\dots,H_{X,n}\}$  
denote the arrangement obtained
from $\A$ by perturbing the lines so that all lines except those
passing through the vertex corresponding to $X$ are in general
position.  

The group of the arrangement $\A_X$ has presentation
\[
G_X=\langle \c_1,\dots ,\c_n \mid  [\c_X,\c_i]=1 
\text{ for } i\in X \text{ and } 
[\c_j,\c_i]=1 \text{ for } j\notin X, i\in [n]\rangle, 
\]
where $\c_{X}=\prod _{i\in X} \c_i$.  
If $|X|=m$, then clearly, $G_X\cong\F_{m-1}\times \Z^{n-m+1}$.   
Let $B_X=B(G_X)$ be the ``local'' Alexander invariant associated to
$X\in L_2$. Note that $B_X$ is trivial if $m=2$.  
For $|X|=m\ge 3$, set 
\[
V_X = \{(t_1,\dots,t_n)\in \T \mid t_X-1=0 \text{ and } 
t_{j}-1=0  \text{ for } j\notin X\}.   
\]
As in the above examples, we have 
\[
V_k(\A_X)=V_X \,\text{ for }\, 1 \le k \le m-2,
\]
$V_k(\A_X)=\bone$ for $m-1 \le k \le \binom{m-1}{2}$, and 
$V_k(\A_X)=\emptyset$ for $k > \binom{m-1}{2}$.

\subsection{Decomposable Alexander Invariants}
\label{subsec:Bcc}
Given an arrangement $\A$ in $\C^2$ that is transverse to infinity, 
let $\A^{\loc} = \prod_{X\in L_2} \A_X$ denote the product 
(see~\cite{OT}) of the 
arrangements $\A_X$ constructed above.
Define the {\em coarse (combinatorial) Alexander invariant} of $\A$
to be the module $\locB(\A) = B(\A^{\loc})\otimes_{\Z\Z^{ns}} \Z\Z^n$ 
induced from the Alexander invariant of $\A^{\loc}$ by
the projection $t_{X,j} \mapsto t_j$.  It is readily seen that  
\[
\locB (\A)=\bigoplus_{X\in L_2} B_X.
\]

In \cite{CSai}, we defined a homomorphism $\Pi:B\to \locB$. 
This map is always surjective, but is not in general a bijection.  
The failure of injectivity is measured by the cokernel of a certain 
linear map $\overline \Psi_3:\overline C_3\to 
\bigoplus_X \overline C_2[X'] \wedge \overline C_1$, where 
$\overline{C}$ denotes the image of a (free) $\L$-module $C$ 
under the augmentation map $\epsilon:\Lambda\to \Z$.  The map 
$\overline \Psi_3$ is determined solely by the intersection lattice, 
see~\cite{CSai}, Section~7.6 for details.  If $\coker\overline \Psi_3=0$, 
and thus $B\cong \locB$, we say that $\A$ is {\em decomposable}. 

\begin{rem}
\label{rem:eisenbud}
In \cite{CSai}, the implications of the surjectivity of the map 
$\overline\Psi_3$ are stated in terms of the $I$-adic completions 
of the modules $B$ and $\locB$.  That these implications apply to 
the modules themselves follows from the ``reflection of isomorphism 
from the completion'' discussed in \cite{Ei},~Exercise~7.5.   
\end{rem}

Now let $V^\loc_k(\A)=V_k(\locB)$ be the {\em $k^{\text{th}}$ coarse 
characteristic variety} of $\A$.  Since the map $\Pi:B\to\locB$ is 
surjective, we have $V^\loc_k(\A) \subseteq V_k(\A)$.  We shall 
refer to the irreducible components of $V^\loc_k(\A)$ as 
{\em local} components of $V_k(\A)$, and to the other irreducible 
components of $V_k(\A)$ as {\em non-local} components.  
If $\A$ is decomposable, then $V^\loc_k(\A) =V_k(\A)$.  
An inductive argument, analogous to the proof of 
Theorem~\ref{thm:cvprod}, yields:
\begin{thm} 
\label{thm:Vsum}  
The coarse characteristic varieties of an arrangement $\A$ are determined 
by the lattice of $\A$, as follows:
\[
V^\loc_k(\A) = \bigcup_{X\in L_2} V_k(G_X) = 
\bigcup_{|X| \ge k+2} V_X.  
\]
In particular, if $\A$ is decomposable, then 
$V_k(\A)=\bigcup_{|X| \ge k+2} V_X$. 
\end{thm}

\begin{rem}
\label{rem:polymat}
The module $B_X=B(G_X)$ depends only on the cardinality $|X|$
of the vertex set $X$.  Consequently, the module $\locB$ depends
only on the number and multiplicities of the elements of $L_2(\A)$.  
On the other hand, the monomial isomorphism type of the coarse 
characteristic variety $V^\loc_k(\A)$ depends on more information 
from the lattice, as the following example of Falk demonstrates. 
\end{rem}

\begin{exm} 
\label{exm:falk}
Consider the two arrangements, $\A_1$ and $\A_2$, 
from Example~4.10 in \cite{Fa}.  
It is readily seen that both arrangements are decomposable:   
$B(\A_{1})=B_{\{1,2,3\}}\oplus B_{\{1,4,5\}}\oplus B_{\{3,5,6\}}
\oplus B_{\{4,6,7\}}$ 
and 
$B(\A_{2})=B_{\{1,2,3\}}\oplus B_{\{1,4,5\}}\oplus B_{\{3,5,6\}}
\oplus B_{\{1,6,7\}}$.   Thus, $B(\A_{1})\cong B(\A_{2})$.  
From Theorem~\ref{thm:Vsum}, we get:
\begin{align*}
V_1(\A_{1})&=V_{\{1,2,3\}}\oplus V_{\{1,4,5\}}\oplus V_{\{3,5,6\}}
\oplus V_{\{4,6,7\}},\\
V_1(\A_{2})&=V_{\{1,2,3\}}\oplus V_{\{1,4,5\}}\oplus V_{\{3,5,6\}}
\oplus V_{\{1,6,7\}}.
\end{align*}
Thus, the varieties $V_{1}(\A_{1})$ and $V_{1}(\A_{2})$ are 
(abstractly) isomorphic.  Nevertheless, there is no monomial 
isomorphism $(\C^{*})^{7}\to (\C^{*})^{7}$ taking  
$V_{1}(\A_{1})$ to $V_{1}(\A_{2})$.  This is proved exactly 
as in \cite{Fa}, using the fact that the ``polymatroids'' 
associated to $\A_{1}$ and $\A_{2}$ are distinct. 
\end{exm}

\section{Cohomology Support Loci}
\label{sec:JumpLoci}

In this section, we use recent results of Arapura~\cite{Ar}
and Libgober \cite{L3} to show that each positive-dimensional, 
irreducible component of the characteristic variety $V_{k}(\A)$ 
is a subtorus of $\T$.  Furthermore, we identify the tangent cone 
of the variety $V_{k}(\A)$ at the point $\bone$.

\subsection{Rank One Local Systems}
\label{sec:LocSys}
Let $\A=\{H_1,\dots ,H_n\}$ be a complex hyperplane arrangement.   
Let $M$ be the complement of $\A$, and let 
$G=\langle \c_{1},\dots , \c_{n} \mid r_{1},\dots , r_{b}\rangle $ 
be the braid presentation for its fundamental group, as in \eqref{eq:bmpres}.  
The generators $\c_i$ are represented by meridional loops around 
$H_i$, with orientations given by the complex structure.  
Since all relations $r_k$ are commutators, $H_1(M; \Z)$ is 
isomorphic to $\Z^n=\langle t_1,\dots ,t_n \mid [t_i,t_j]=1\rangle$, 
with identification given by $\c_i\mapsto t_i$.  
Furthermore, $H^1(M; \C^*)$ is isomorphic to the 
algebraic torus $\T$, with coordinates $\bt=(t_1,\dots ,t_n)$.   

Each point $\bt$ in $\T$ determines a representation
$G\to\C^*=\GL(1,\C)$, $\c_i\mapsto t_i$, and an associated rank 
one local system, which we denote by $\C_\bt$.  
For generic $\bt$, the (co)homology of $M$ with coefficients in
$\C_\bt$ vanishes.  Those $\bt$ for which $H^{r}(M;\C_\bt)$ does not 
vanish comprise the {\em cohomology support loci}
\[
W^{r}_k(M)=\{\bt\in\T \mid \rank H^{r}(M; \C_{\bt})\ge k\}.
\]
These loci are algebraic subvarieties of $\T$, which depend 
only on the homotopy type of $M$, and a generating set for $G=\pi_1(M)$.  
We now relate $W^1_k(M)$ to the characteristic variety $V_k(\A)$.

If $\A$ is a {\em general position} (through rank~$2$) arrangement, 
then the group of $\A$ is free abelian and the Alexander invariant 
is trivial.  Thus, $V_{k}(\A)=\emptyset$ in this instance.  
On the other hand, $W^1_{k}(M)=\{\bone\}$ for $1 \le k \le n$ and 
$W^1_{k}(M)=\emptyset$ for $k>n$ (see e.g.~\cite{Hat}).

If $\A$ is not a general position arrangement, 
we have the following.  As shown by Hironaka \cite{Hir} (see also \cite{L3}), 
an analogous result holds for an arbitrary CW-complex $M$ with 
torsion-free first homology (and non-trivial Alexander invariant).

\begin{thm} 
\label{thm:vw}
Let $\A$ be an affine arrangement of $n$ lines in $\C^2$, with complement $M$.          
Set $N=\min\{ n, \binom{n}{2} - b_2(M) \}$.   If $H_1(M')\ne 0$, 
then, for $1 \le k \le N$, the characteristic variety $V_k(\A)$ 
coincides with the cohomology support locus $W^1_k(M)$.
\end{thm}
\begin{proof}
For an arbitrary arrangement, the homology of $M$ with
coefficients in the local system $\C_\bt$ is naturally isomorphic 
to the homology of the chain complex 
$C_{\*}(M;\C_\bt):=C_\*(M') \otimes_{\L} \C_\bt$, where 
$\L=\Z\Z^{n}$ and $(C_{\*}(M'),\partial_{\*})$ is the chain complex
of the maximal abelian cover of $M$ specified in \eqref{eq:m'cx}.  
The terms of $C_{\*}(M;\C_\bt)$ are finite-dimensional complex vector spaces, 
and the boundary maps, $\bp_k(\bt)$, are the evaluations of $\partial_k$ 
at $\bt$.  

Let $W_1^k(M)=\{\bt \in \T \mid \rank H_1(M;\C_\bt)\ge k\}$ denote the 
homology support loci of $M$.  Evidently, these varieties coincide with the 
cohomology support loci defined above.  From the above description of the 
complex $C_{\*}(M;\C_\bt)$, it is clear that $W_1^n(M)=\{\bone\}$ 
and $W_k(M)=\emptyset$ for $k>n$, while 
as noted in \eqref{eqn:vbigk}, $V_k(\A)=\{\bone\}$ for 
$n \le k \le \binom{n}{2}-b$, where $b=b_2(M)$.  Thus if
$n\le \binom{n}{2}-b$, we have $V_n(\A)=W_1^n(M)$.

For $k<n$, since $\bp_1(\bt)=\mathbf{d}_1(\bt)$ is of rank at most one, 
the homology support loci are given by 
$W_1^k(M)=\{\bt \in \T \mid \rank \bp_2(\bt) \le n-k-1\}$.  
Recall that the Alexander invariant $B=B(\A)$ is presented by the matrix 
$\Delta$ of Theorem~\ref{thm:aipresthm}.  Upon tensoring the terms of 
the commuting diagram~\eqref{eq:mapconediagram} with $\C$, and evaluating 
the maps at $\bt$, it is readily seen that 
\[
W_1^k(M)=\{\bt \mid \rank \bp_2(\bt) \le n-k-1\}=
\{\bt \mid \rank\bD(\bt) \le \tbinom{n}{2}-k\}=V_k(\A)
\]
for $1\le k \le N$.
\end{proof}

We now pursue a more precise description of the varieties $V_k(\A)$.
First, note that $M=\CP^{2}\setminus \bigcup_{H\in\A^*} H$, where $\A^*$ 
is the projectivization of the cone $\cA$.  By blowing up the 
singularities, we see that $M$ is biholomorphically equivalent 
to the complement of a normal-crossing divisor in a smooth, 
simply-connected projective variety.  
(This holds more generally for complex subspace arrangements, 
as shown by De~Concini and Procesi~\cite{DCP}.)   
The following result of Arapura~\cite{Ar} describes 
the support loci of such quasiprojective varieties.  

\begin{thm}[Arapura~\cite{Ar}] 
\label{thm:donu}  
Let $M$ be the complement of a normal-crossing divisor in a compact 
K\"{a}hler manifold with vanishing first homology.  The cohomology 
support locus $W^1_{k}(M)$ is a finite union of torsion-translated 
subtori of the algebraic torus $\T$, where $n=b_{1}(M)$.
\end{thm}

In other words, each irreducible component of $W^1_{k}(M)$ is of 
the form $\mathbf{q}\cdot (\C^*)^{r}$, for some integer $r\ge 0$ 
and some torsion element $\mathbf{q}=(q_1,\dots,q_n)\in\T$.  
In \cite{L2} and \cite{L3}, Libgober studies in detail the case 
where $M$ is the complement of a plane algebraic curve $\CC$.  
He shows that, if $r>0$, the order of each root of unity $q_{j}$ divides 
the greatest common divisor of the degrees of the branches of $\CC$.  
In particular, if $\CC$ is a line arrangement, then $\mathbf{q}=\bone$.  
As a result, we have the following.  

\begin{thm}[Libgober~\cite{L3}] 
\label{thm:anatoly} 
Let $\A$ be an arrangement of $n$ complex hyperplanes.  
Then each positive-dimensional component of the 
characteristic variety $V_{k}(\A)$ is a 
connected 
subtorus 
of the algebraic torus $\T$.
\end{thm}

In view of this theorem, define the {\em $k^{\text th}$ central
characteristic subvariety} of the arrangement $\A$ to be the subvariety 
$\check{V}_k(\A)$ of $\T$ consisting of all irreducible 
components of $V_k(\A)$ passing through $\bone$.  Note that 
$V_{k}(\A)\setminus \check{V}_{k}(\A)$ consists of a finite set 
of nontrivial torsion elements.  In all examples we have examined, 
$\check{V}_1(\A)=V_1(\A)$.  In general however, 
$\check{V}_{k}(\A)\subsetneqq V_{k}(\A)$, as the following example 
illustrates.  

\begin{exm}
\label{exm:diamond} 
Consider the central $3$-arrangement $\A=\{H_1,\dots ,H_7\}$, 
with defining polynomial $Q=x(x+y+z)(x+y-z)y(x-y-z)(x-y+z)z$.  
The lattice of $\A$ has six rank two flats of multiplicity three.   
The first characteristic variety, $V_1(\A)$, is the union of 
nine $2$-dimensional subtori of $(\C^*)^7$, six of which are local 
(see Theorem~\ref{thm:Vsum}).  The remaining three $2$-tori, 
\begin{align*}\label{nonlocdiam}
&\{\bt\in (\C^*)^7 \mid t_1=t_4,\, t_2=t_3,\, t_5=t_7,\, t_6=1,\, t_1t_2t_5=1\}, \\
&\{\bt\in (\C^*)^7 \mid t_1=t_5,\, t_2=t_6,\, t_4=t_7,\, t_3=1,\, t_1t_2t_4=1\}, \\
&\{\bt\in (\C^*)^7 \mid t_1=t_7,\, t_3=t_6,\, t_4=t_5,\, t_2=1,\, t_1t_3t_4=1\},
\end{align*}
arise from subarrangements that are isomorphic to the braid arrangement $\A_4$ 
(see Sections~\ref{sec:essentialtori} and~\ref{sec:braidarrangement} below).  
The second characteristic variety, $V_2(\A)$, consists of the two points
of intersection of the three $2$-tori above, $\bone$  and  $(-1,1,1,-1,-1,1,-1)$,  
whereas $\check{V}_{2}(\A)=\bone$.  
\end{exm}

\subsection{Tangent Cones}
\label{sec:TanCone}
We now study the tangent cone at the origin $\bone\in\T$ of the 
variety $V_k(\A)=W^1_k(M)$.  The tangent space of $H^1(M,\C^*)=\T$ 
at $\bone$ is $H^1(M,\C)=\C^n$, with coordinates 
$\bl=(\l_1,\dots ,\l_n)$. The exponential map 
$\TT_{\bone}\left( \T \right)\to \T$ 
is just the coefficient map $H^1(M,\C)\to H^1(M,\C^*)$ induced 
by $\exp: \C\to\C^*$, $\l_i \mapsto e^{\l_i}=t_i$.  

Let $\cV_{k}(\A)$ denote the tangent cone at $\bone$ of the variety 
$V_{k}(\A)$.  Clearly, this coincides with the tangent cone at $\bone$ 
of $\check{V}_{k}(\A)$.  
From Theorem~\ref{thm:anatoly}, we know that $\check{V}_k(\A)$ 
is a (central) arrangement of subtori in $\T$.  Thus, $\cV_k(\A)$ 
is a central arrangement of subspaces in $\C^n$.  We shall find 
defining equations for this variety, and relate them 
to the presentation of the Alexander invariant 
of $M$ from Theorem~\ref{thm:aipresthm}. 

Let $f\in \L$ be a Laurent polynomial, and $\mathbf{f}:\T\to\C$ 
the corresponding map, defined in Section~\ref{subsec:cva}. 
The derivative of this map at the identity, 
$\mathbf{f}_*:\TT_{\bone}\left(\T\right)\to\C$, is given by 
$\mathbf{f}_*(\bl) = {\frac{d}{dx}}\bigr\rvert _{x=0} 
f(e^{x\l_1},\dots,e^{x\l_n}).$
If $f$ and $g$ are two Laurent polynomials, the Product Rule yields  
$(\mathbf{fg})_*(\bl)=\mathbf{f}_*(\bl)\mathbf{g}(\bone)+
\mathbf{f}(\bone)\mathbf{g}_*(\bl).$ 
More generally, for $F \in \Mat_{p\times q}(\L)$ 
and $G \in \Mat_{q\times r}(\L)$, matrix multiplication 
and the differentiation rules yield
\begin{equation}  
\label{eq:leibniz}
(\mathbf{F\cdot G})_*(\bl)=\mathbf{F}_*(\bl)\cdot\mathbf{G}(\bone)+
\mathbf{F}(\bone)\cdot\mathbf{G}_*(\bl).
\end{equation}

Recall the boundary map $d_k:C_k \to C_{k-1}$ 
from the standard resolution \eqref{eq:znres}.  Let $\mathbf{d}_{k}$ 
be the corresponding map, and $\bd_k=(\mathbf{d}_k)_{*}$ its 
derivative at $\bone$.   Note that 
$\mathbf{d}_k(\bone)=\bz$ is the zero matrix, and that 
$\bd_k(\l_1,\dots ,\l_n)=\mathbf{d}_k(1-\l_1,\dots , 1-\l_n)$.  

Recall also the presentation matrix, 
$\Delta = \left(\smallmatrix\Phi \\ d_3 
\endsmallmatrix\right):\L^p\to \L^q$, 
for the Alexander invariant of $M$ from Theorem~\ref{thm:aipresthm}, 
with $p=b_2(M)+\binom{n}{3}$ and $q=\binom{n}{2}$.  For
$\bt \in \T$,  let $\bD(\bt):\C^p\to\C^q$ be the 
corresponding map, given by
\begin{equation}  
\label{eq:delt}
\bD(\bt) = \begin{pmatrix}\boldsymbol{\Phi}(\bt) \\
\mathbf{d}_3(\bt) \end{pmatrix}.
\end{equation}
As noted above, we have $V_k(\A) = 
\{\bt \in \T \mid \rank \bD(\bt) \le \binom{n}{2}-k\}$.
We now obtain an analogous description of the tangent cone 
$\cV_k(\A)$. 

\begin{thm} 
\label{thm:tj}
Let $\A$ be an arrangement of $n$ hyperplanes, with complement $M$, 
and non-trivial Alexander invariant $B$, presented by 
$\Delta = \left(\smallmatrix \Phi \\ d_3 \endsmallmatrix\right)$. 
Then the tangent cone, $\cV_k(\A)$, of the characteristic variety $V_k(\A)$ 
at the point $\bone$ is a subspace arrangement in $\C^n$, defined 
by the equations
\begin{equation*}  
\label{eq:tjA}
\cV_{k}(\A)=
\biggl\lbrace
\bl\in\C^n \ \bigl\vert\ 
\rank \begin{pmatrix}\boldsymbol{\Phi}(\bone) \\
\bd_3(\bl) \end{pmatrix} \le \binom{n}{2}-k\biggr.
\biggr\rbrace.
\end{equation*}
\end{thm}
\begin{proof}
The identification in Theorem~\ref{thm:vw} of the characteristic variety 
$V_k(\A)$ and the (co)homology support locus shows that 
$V_k(\A) =\{\bt\in\T \vert \rank \bp_2(\bt) \le n-k-1\}$, 
where 
$\bp_2(\bt)$ is obtained from the Alexander matrix $\p_2$ 
of~\eqref{eq:m'cx} by evaluation at $\bt$.  Thus the tangent 
cone at $\bone$ is given by 
$\cV_k(\A)=\{\bl\in\C^n \mid \rank (\bp_2)_*(\bl) \le n-k-1\}$.

Recall that the image of a free $\L$-module $C$ under the 
augmentation map is denoted by $\overline{C}$.  
Also, recall from~\eqref{eq:phid2} that $\p_2 = d_2 \circ \Phi$. 
Consequently, for each $\bt$, we have 
$\bp_2(\bt) = \boldsymbol{\Phi}(\bt) \cdot \mathbf{d}_{2}(\bt)$.  
Using~\eqref{eq:leibniz} and the fact that $\mathbf{d}_2(\bone)=\bz$, 
we obtain $(\bp_2)_*(\bl) = \boldsymbol{\Phi}(\bone) \cdot \bd_2(\bl)$.
Thus, we have a commuting diagram
\begin{equation*}
\begin{CD}
&& \displaystyle{\bigoplus_{X\in L_2}} \overline{C}_{1}[X'] @>(\bp_{2})_*(\bl)>> 
\overline{C}_{1} @>\bd_{1}(\bl)>> \overline{C}_{0}\\
&& @VV{\boldsymbol{\Phi}(\bone)}V  @VV{\id}V  @VV{\id}V \\
\overline{C}_{3} @>\bd_{3}(\bl)>> \overline{C}_{2} @>\bd_{2}(\bl)>>
\overline{C}_{1}@>\bd_{1}(\bl)>>\overline{C}_{0}\\
\end{CD}
\end{equation*}
and an exercise in homological algebra, using a mapping cone
construction as in the proof of Theorem~\ref{thm:aipresthm}, 
completes the proof.
\end{proof}

\begin{rem} 
\label{rem:phiat1}  
Recall from~\ref{subsec:AlexInv} that the map 
$\Phi:\bigoplus_{X \in L_2} C_1[X'] \to C_2$ is defined by 
$\Phi|_{C_1[X']}=
\Theta_{2}(\d)\circ\Phi(A_X)\circ\Theta_{1}(\d^{-1})|_{C_1[X']}$, 
where $A_X^\d$ is the braid monodromy generator corresponding 
to $X\in L_2$.  Since pure braids are IA-automorphisms of the 
free group, we have $\boldsymbol{\Theta}_k(\d)(\bone)=\id$ 
for all $\d\in P_n$.  Thus, 
$\boldsymbol{\Phi}(\bone)|_{\overline{C}_1[X']}=
\boldsymbol{\Phi}_X(\bone)$, and it is readily checked that 
the latter is given by
$\boldsymbol{\Phi}_X(\bone)(e_i)=e_i \wedge \overline\nabla_X$ for 
each $i\in X'$, where $\overline\nabla_X = \sum_{j\in X} e_j$.
It follows that the map $\boldsymbol{\Phi}(\bone)$ is determined by
the vertex sets $X\in L_2$.  However, it is not a priori clear that 
this map is combinatorial, 
as the vertex sets and their ordering are dictated by the braid monodromy. 
By comparison with the Orlik-Solomon algebra, the combinatorial
nature of the map $\boldsymbol{\Phi}(\bone)$ will be established in
the next section.  
\end{rem}

\section{Relation to the Orlik-Solomon Algebra}
\label{sec:OSalgebra}

In this section, we compare Falk's invariant of the Orlik-Solomon 
algebra of an arrangement of hyperplanes $\A$ with the characteristic 
varieties of the group of $\A$.

\subsection{Orlik-Solomon Algebra}
\label{sec:O-S}
Let $\A=\{H_1,\dots,H_n\}$ be a (central) arrangement.  
Let $E$ be the graded exterior algebra over $\C$ generated 
by $1$ and $e_1,\dots,e_n$.  Note that $E^r = \overline{C}_r$.  
Define $\overline{d}_{r}:E^r\to E^{r-1}$ by
$\overline{d}_{r}(e_{J})=\sum _{q=1}^{r} (-1)^{r} e_{J\setminus\{j_q\}}$.  
Let $I=I(\A)$ be the homogeneous ideal in $E$ generated by 
$\{\overline{d}(e_{J}) \mid \dim \bigcap_{j\in J} H_{j} <|J|\}$.
The Orlik-Solomon algebra of $\A$ is the graded algebra $A=E/I$.

Let $p:E \to A$ denote the natural projection, and write 
$a_{i}=p(e_{i})$.  Note that $I^{0}=I^{1}=0$, and so the maps 
$p:E^r\to A^r$ are isomorphisms for $r\le 1$.  
An {\bf nbc} (no broken circuit) basis for $A^2$, corresponding 
to the ordering $H_1,\dots,H_n$ of the hyperplanes, consists of 
all generators $a_j\wedge a_k$, except those for which there 
exists $i<\min\{j,k\}$ such that $H_i\cap H_j \cap
H_k \in L_2(\A)$, see \cite{OT}.  With this choice of basis, 
for $H_i\cap H_j \cap H_k \in L_2(\A)$,  
the projection $p:E^2\to A^2$ is given by 
$p(e_i\wedge e_j)=a_i\wedge a_j$, 
$p(e_i\wedge e_k)=a_i\wedge a_k$, and
$p(e_j\wedge e_k)=a_i\wedge a_k-a_i\wedge a_j$.  

Associated to each $\bl \in \C^n$, we have an element 
$\omega = \sum_{i=1}^n \l_i a_i$ of $A^1$.  Left-multiplication by 
$\omega$ defines a map $\mu:A^{r}\to A^{r+1}$.  Clearly, 
$\mu\circ\mu=0$, so $(A,\mu)$ is a complex.  
The cohomology support loci of the OS-algebra are Falk's invariant
varieties $\RR^{r}(A)=\{\bl\in\C^n\mid H^{r}(A,\mu)\neq 0\}$
(denoted in \cite{Fa} by $\RR_{r}(A)$).  We shall consider here
only the first such locus, filtered by the family of subvarieties
$\RR^{1}_k(A)=\{\bl\in\C^n\mid \rank H^1(A,\mu)\ge k\}$.  
See \cite{OT} and \cite{Fa} for detailed discussions of 
the OS-algebra and the variety $\RR^1(A)$.

If the arrangement $\A$ is in general position (through rank 2), 
it is well-known that $\RR^1(A)=\{\bz\}$, see e.g.~\cite{Yuz}.  
If $\A$ is not in general position, we have the following.

\begin{thm} 
\label{thm:randv}
Let $\A$ be a central arrangement of complex hyperplanes with
non-trivial Alexander invariant.  Then the cohomology support locus
$\RR^{1}_k(A)$ of the Orlik-Solomon algebra of $\A$ coincides with the
tangent cone $\cV_k(\A)$ of the characteristic variety $V_k(\A)$ 
at the point $\bone$.
\end{thm}
\begin{proof}
Given $\bl\in\C^n$, let $\hat\omega=\sum\l_i e_i \in E^1$.  
Left-multiplication by 
$\hat\omega$ defines a map $\hat\mu:E^{r}\to E^{r+1}$, and as above,
$(E,\hat\mu)$ is a complex.  The natural projection $p:E\to A$ 
is a chain map with kernel $I$. If $\bl\neq \bz$, the complex 
$(E,\hat\mu)$ is acyclic, and we have $H^r(A,\mu)=H^{r+1}(I,\hat\mu)$.  
Thus to study $H^*(A,\mu)$, $\RR^{1}_{k}(A)$, etc., one can pass from 
the OS-algebra to the OS-ideal, see~\cite{Yuz}, \cite{Fa}. 
However, to facilitate comparison with the Alexander invariant, we
opt for a slightly different approach.

Consider the mapping cone $P$ of $p$.  This is a complex with terms 
$P^{r}=A^{r-1} \oplus E^{r}$, and differentials $\xi:P^{r}\to P^{r+1}$ 
given by $\xi(x,y) = (p(y)-\omega\wedge x,\hat\omega\wedge y)$. The
natural inclusion and projection yield an exact sequence of complexes
$0\to A^{\*-1}\to P^\* \to E^\*\to 0$.  As above, we have
$H^r(A,\mu)=H^{r+1}(P,\xi)$ if $\bl\neq \bz$.  Thus we can write
$\widehat\RR^{1}_k(A)=\{\bl\in\C^n\setminus\{\bz\}\mid \rank H^2(P,\xi)\ge k\}$,
where $\widehat\RR^{1}_k(A)=\RR^{1}_k(A)\setminus\{\bz\}$.

Define an automorphism
$\psi:P^2\to P^2$ by $\psi(x,y)=(x,y-\hat\omega\wedge x)$.
Then, using the fact that $A^r=E^r$ for $r\le 1$,
we see that $\psi\circ\xi:P^1\to P^2$ is given by
$(x,y)\mapsto (y-\omega\wedge x,0)$, and
$\xi\circ\psi^{-1}:P^2\to P^3$ by
$(x,y)\mapsto (p(y),\hat\omega\wedge y)$.
It follows that $H^2(P,\xi) = \ker\phi$, where
$\phi=(p,\hat\mu):E^2\to A^2\oplus E^3$ is the restriction of
$\xi\circ\psi^{-1}$ to $E^2\subseteq P^2$.  Thus $\widehat\RR^{1}_k(A)$ consists
of all $\bl\neq \bz$ for which $\dim\ker\phi\ge k$, i.e., for which
$\rank\phi(\bl) \le \binom{n}{2}-k$.  For an arrangement $\A$ 
with non-trivial Alexander invariant, it is readily checked that
$\rank\phi(\bz) < \binom{n}{2}$.  For such $\A$, we have
$\RR^{1}_k(A) = \{\bl \in \C^n \mid \rank \phi(\bl) \le \binom{n}{2}-k\}$.

The choice of {\bf nbc} basis for $A^2$ above provides
a natural identification $A^2 \leftrightarrow
\bigoplus_{X\in L_2} \overline{C}_1[X']$.  Using this
identification, the fact that $E^k=\overline{C}_k$,
and the description of the map
$\boldsymbol{\Phi}(\bone)$ from Remark~\ref{rem:phiat1},
it is readily checked that $p$ is dual to $\Phi(\bone)$
and $\hat\mu(\bl)$ is dual to $\bd_3(\bl)$.  It follows that
$\phi(\bl)=(p,\hat\mu(\bl))$ is the transpose of the matrix  
$\Bigl(\smallmatrix
\boldsymbol{\Phi}(\bone)\\\bd_3(\bl) \endsmallmatrix\Bigr)$
from Theorem~\ref{thm:tj}.  Thus, $\RR^{1}_k(A)=\cV_k(\A)$.
\end{proof}

Recall from Theorem~\ref{thm:anatoly} that all irreducible components 
of the variety $V_{k}(\A)$ passing through $\bone$ are subtori of $\T$.  
Thus, all irreducible components of the tangent cone $\cV_{k}(\A)$ 
are (linear) subspaces of $\C^{n}$.  As a consequence of 
Theorem~\ref{thm:randv}, we have the following, which 
resolves a conjecture made by Falk (\cite{Fa}, Conjecture~4.7).
\begin{cor} 
\label{cor:falkconj}
The variety $\RR^{1}_{k}(A)$ is the union of a subspace arrangement.
\end{cor}

Since the subspace arrangement $\cV_{k}(\A)=\RR^{1}_{k}(A)$ is determined, 
up to a linear change of basis in $\C^{n}$, by the intersection lattice 
$L(\A)$, the central arrangement of tori $\check{V}_{k}(\A)$ is determined, 
up to a monomial change of basis in $\T$, by $L(\A)$.  Thus, as another 
consequence of Theorem~\ref{thm:randv}, we have the following.
\begin{cor} 
\label{cor:vcomb}
The monomial isomorphism type of the central characteristic subvariety 
$\check{V}_{k}(\A)$ is a combinatorial invariant of the arrangement $\A$.
\end{cor}

\begin{rem} 
Libgober proves Corollaries~\ref{cor:falkconj} and \ref{cor:vcomb} by 
other means in \cite{L3}.  However, the assertion there that the (entire) 
characteristic variety $V_k(\A)$ is combinatorial does not follow immediately 
from these results.  As Example~\ref{exm:diamond} illustrates, the variety 
$V_k(\A)$ may contain isolated points, which are not a priori 
combinatorially determined. 
\end{rem}

\subsection{Essential Tori and Non-Local Components}
\label{sec:essentialtori}
We use the identification of the variety $\RR^{1}(A)$ and the 
tangent cone of the first characteristic variety provided by 
Theorem~\ref{thm:randv} to obtain a combinatorial means for detecting 
components of $\check{V}_1(\A)$.

Positive-dimensional, irreducible components of the 
characteristic variety $V_{k}(\A)$ which are 
not contained in any coordinate torus $t_i=1$ of $\T$ are called 
{\em essential tori} in \cite{L3}.  For an arrangement $\A$ of 
cardinality $n$, with rank two lattice elements of multiplicities 
only $m$ and $2$, Libgober gives an algorithm for detecting such 
essential tori, see~\cite{L3}, Section 3.  Ingredients of this 
algorithm include polytopes of quasi-adjunction, and the 
calculation of the cohomology of a certain ideal sheaf.  

Using Theorem~\ref{thm:randv} above, we obtain a combinatorial 
alternative to this approach.  While this alternative method applies 
only to the first characteristic variety $V_1(\A)$, it has the advantage of being 
applicable to an arbitrary arrangement.  To describe this method, we 
must first recall Falk's description of the variety $\RR^1(A)$ from~\cite{Fa}.

Let $\A=\{H_1,\dots,H_n\}$ be an arrangement with lattice $L=L(\A)$, 
and consider $\C^n$ with basis $e_i$ and coordinates $\lambda_i$ 
corresponding to the elements of $\A$ (resp.~of $L_1$).  
A partition $\Pi$ of $[n]$ (or of $L_1$) is 
{\em neighborly} if, for all flats $X\in L_2$ and all blocks $\pi$ of 
$\Pi$, $|\pi \cap X| \ge |X|-1$ implies that $X \subseteq \pi$.  
Partitions with more than one block will be called non-trivial.  
A flat $X$ that is contained in a single block of $\Pi$ is said 
to be {\em monochrome}, and is called {\em polychrome} otherwise.  
Note that all rank two flats of multiplicity two are necessarily 
monochrome.

For each flat $X\in L_2$, let $H_X$ be the hyperplane in $\C^n$ 
defined by 
$\sum_{i\in X} \lambda_i=0$, and let $H_{[n]}=\{\sum_{i=1}^n 
\lambda_i=0\}$.  
If $\Pi$ is a neighborly partition, let $S_\Pi$ denote the subspace 
of $\C^n$ defined by $S_\Pi=H_{[n]} \cap \bigcap H_X$, where $X$ 
ranges over all polychrome flats in $L_2$.  Associated to $\Pi$,
we also have a vector-valued skew-symmetric bilinear form 
$\langle\bl ,\bl' \rangle_\Pi$ on $\C^n$, the components of which 
are the $2\times 2$ determinants 
$\left|\smallmatrix \lambda_i&\lambda_j\\ 
\lambda'_i&\lambda'_j\endsmallmatrix \right|$, for 
$i,j\in\pi \subset \Pi$.  Falk's characterization of the 
variety $\RR^1(A)$ is the following.

\begin{thm}[Falk~\cite{Fa}, Theorem~3.10] 
\label{thm:falk}
$\bl \in \RR^1(A)$ if and only if there exists a subarrangement $\A'$ 
of $\A$ and a neighborly partition $\Pi$ of $L_1(\A')$ such that
\begin{romenum}
\item $\bl \in S_\Pi$; and
\item there exists $\bl' \in S_\Pi$, not proportional to $\bl$, 
such that $\langle \bl,\bl'\rangle_\Pi=\mathbf{0}$.
\end{romenum}
\end{thm}

If $\Pi$ is a neighborly partition, let $\cV_\Pi$ denote the 
subvariety of the subspace $S_\Pi$ on which the form 
$\langle \ ,\ \rangle_\Pi$ is degenerate.  In all known examples, 
$\cV_\Pi = S_\Pi$, see~\cite{Fa} and below.  By Theorem~\ref{thm:randv}, 
the variety $\cV_\Pi$ is linear in general.  This does not, however, 
rule out the possibility that the containment $\cV_\Pi \subset S_\Pi$ 
is strict.

\begin{exm}\label{exm:hessian} 
Let $\A$ be the Hessian configuration, with defining polynomial
$Q=x_1x_2x_3\prod_{i,j=0,1,2} (x_1+\omega^i x_2+\omega^j x_3)$, where 
$\omega=\exp(2\pi\ii/3)$.  The characteristic varieties $V_k(\A)$ were determined
by Libgober in~\cite{L3}, Example~5.  In particular, he finds an 
essential torus $V \subset V_1(\A)$ of dimension three.  Consequently, the 
variety $\RR^1(A)$ also has a three-dimensional (linear) component $S$.  
This component is readily recovered using Theorem~\ref{thm:falk}. 
Note that $L_2(\A)$ consists of $9$~flats of multiplicity $4$ and $12$ 
flats of multiplicity $2$.  
Let $H_i=\ker x_i$ and let 
$H_{i,j}=\ker(x_1+\omega^i x_2+\omega^j x_3)$.  The partition 
$\Pi=\big( H_1,H_2,H_3 \ \big|\ H_{0,0},H_{1,2},H_{2,1} \ \big|\ H_{0,1},
H_{1,0},H_{2,2} \ \big|\ H_{0,2},H_{1,1},H_{2,0}\big)$ 
is neighborly, and all flats of multiplicity $4$ are polychrome.  
The subspace $S_\Pi$ is three-dimensional.  Checking that the form 
$\langle \ ,\ \rangle_\Pi$ is trivial on $S_\Pi$, we have $S=S_\Pi$.
\end{exm}

\begin{rem}
\label{rem:tutte}
As illustrated by the above example, for an arbitrary arrangement $\A$, 
the variety $\RR^1(A)$ may have non-local components of dimension greater 
than two.  This leaves open the possibility that the Tutte polynomial of 
(the matroid of) an arrangement is not determined by the 
Orlik-Solomon algebra, see~\cite{Fa}, Section~3.
\end{rem}

In practice (see below), we use Falk's characterization of $\RR^1(A)$ 
to detect essential tori and non-local components of $V_1(\A)$.  
Theorems~\ref{thm:randv} and \ref{thm:falk} yield:

\begin{prop}
\label{prop:neighborly}
Let $\A$ be an arrangement of $n$ complex hyperplanes.  If $\Pi$ is an 
associated neighborly partition of $[n]$, and $\cV_\Pi$ is 
the subspace of $\C^n$ on which the form $\langle \ ,\ \rangle_\Pi$ 
is degenerate, then $V_\Pi = \exp(\cV_\Pi)$ is an essential torus of 
the first characteristic variety $V_1(\A)$.
\end{prop}

\section{Monomial Arrangements}
\label{sec:monomial}

In this section, we use the results on characteristic varieties obtained
above to study the arrangements associated with the full monomial 
groups $G(r,1,\ll)$ and the irreducible subgroups $G(r,r,\ll)$, $r\ge 2$. 
Defining polynomials for the reflection arrangements $\A_{r,1,\ll}$ and
$\A_{r,r,\ll}$ corresponding to these groups are given by
\[
Q(\A_{r,1,\ll}) = x_1\cdots x_\ll \prod_{1\le i < j \le \ll}
(x_i^r-x_j^r)
\quad \text{and} \quad
Q(\A_{r,r,\ll}) = \prod_{1\le i < j \le \ll} (x_i^r-x_j^r).
\]
Note that the arrangements $\A_{2,1,\ll}$ and $\A_{2,2,\ll}$ 
are the Coxeter arrangements of type $\Bb$ and $\Dd$ respectively.
The fundamental group, $P(r,k,\ll) = \pi_1(M(\A_{r,k,\ll}))$, is the 
generalized pure braid group associated to the complex reflection 
group $G(r,k,\ll)$, $k=1$, $r$.  
Presentations for the (full) braid 
groups associated to these reflection groups were recently 
found in~\cite{BMR}. 

The arrangement $\A_{r,1,\ll}$ is 
fiber-type, with exponents $\{n_1,n_2,\dots,n_\ll\}$, 
where $n_{i}=(i-1)r+1$. Hence, the group $P(r,1,\ll)$ admits 
the structure of an iterated semidirect product of free groups: 
$P(r,1,\ll) \cong \F_{n_\ll} \rtimes \cdots \rtimes \F_{n_2}\rtimes
\F_{n_1}$.  

The groups $P(r,r,\ll)$ also admit such structure.  Define 
$\pi:\C^{\ll} \to \C^{\ll-1}$ by 
$\pi(x_{1},\dots,x_{\ll-1},x_{\ll}) = 
(x_{1}^{r}-x_{\ll}^{r},\dots, x_{\ll-1}^{r}-x_{\ll}^{r})$.  
The restriction of $\pi$ to the complement $M(\A_{r,r,\ll})$ is a 
generalized Brieskorn bundle (see \cite{Br} for the case $r=2$).  
The base space of this bundle is homotopy equivalent to the complement 
of the braid arrangement $\A_{\ll}$, and the fiber 
is a surface of genus $\frac{r^{\ll-1}(\ll-3)-(r-2)(\ll-1)+2}{2}$ 
with $r^{\ll -1}$ punctures.  It follows that
$P(r,r,\ll) \cong \F_{m_{\ll}} \rtimes \cdots \rtimes 
\F_{m_{2}}\rtimes \F_{m_{1}}$, 
where $m_{i} = i$ for $1 \le i \le \ll-1$, and 
$m_{\ll}=r^{\ll-1}(\ll-2)+(r-2)(\ll-1)+1$.

Let $\Pi(r,1,\ll) = \F_{n_{\ll}} \times \cdots \times \F_{n_{2}}\times
\F_{n_{1}}$ and $\Pi(r,r,\ll) = \F_{m_{\ll}} \times \cdots \times 
\F_{m_{2}}\times \F_{m_{1}}$ denote the corresponding direct products 
of free groups.  

\begin{thm} 
\label{thm:PandPi}
For $\ll \ge 3$, the groups $P(r,k,\ll)$ and $\Pi(r,k,\ll)$ 
are not isomorphic.
\end{thm}

We establish this result by distinguishing the characteristic 
varieties $V_1(P(r,k,\ll))$ and $V_1(\Pi(r,k,\ll))$. 
The latter may be determined using 
Examples~\ref{exm:decone}~and~\ref{exm:Vfreeprod}.  
In particular, this variety has $\ll-1$ irreducible components.  
As we shall see in Proposition~\ref{prop:topmono2} and 
Remark~\ref{rem:fullmono}, the variety $V_1(P(r,k,\ll))$ 
has many more components.  

An analogous result for Artin's pure braid group $P_{\ll+1}$ 
and the direct product 
$\Pi_{\ell+1}=\F_\ll \times \cdots \times \F_2 \times \F_1$ 
was proved in \cite{CSpn} by other means.  
Proposition~\ref{prop:v1braid} below provides another 
proof of this fact.

\subsection{Monomial Arrangements in $\C^3$} 
\label{sec:rr3}
We now determine the structure of the characteristic variety
$V_1(\A)$ of every monomial arrangement $\A=\A_{r,r,3}$ 
in $\C^3$.  We first show that  
the variety $V_1(\A)$ contains 
an essential torus of dimension two.  The cases $r=2$ and $r=3$ 
were considered previously by Falk~\cite{Fa} and Libgober~\cite{L3}.  

A defining polynomial for $\A$ is given by 
$Q(\A)=(x_1^r-x_2^r)(x_1^r-x_3^r)(x_2^r-x_3^r)$.  
Denote the hyperplanes of $\A$ by 
$H_{i,j}^{(k)} = \ker (x_i - \zeta^k x_j)$, 
where $\zeta = \exp(2\pi\ii/r)$, $1\le i < j \le 3$, 
and $1\le k \le r$.  In this notation, the rank two elements 
of the lattice of $\A$ are given by
\[
L_2 = \biggl\lbrace
\bigcap_{k=1}^r H_{i,j}^{(k)}, 1\le i < j \le 3,\quad
H_{1,2}^{(p)} \cap H_{2,3}^{(q)} \cap H_{1,3}^{(s)},
\ \text{where} \ s \equiv p+q \bmod r 
\biggr\rbrace.
\]

Consider the partition $\Pi = \Big( 
H_{1,2}^{(1)},\dots,H_{1,2}^{(r)} \ \Big| \ 
H_{1,3}^{(1)},\dots,H_{1,3}^{(r)} \ \Big| \ 
H_{2,3}^{(1)},\dots,H_{2,3}^{(r)} \Big)$ 
of the rank one elements of $L(\A)$.  
It is readily checked that this partition is neighborly, with 
polychrome flats all of the triple points 
$H_{1,2}^{(p)} \cap H_{2,3}^{(q)} \cap H_{1,3}^{(s)}$ noted above.  
Let $n=3r=|\A|$, and consider $\C^n$ with basis $e_{i,j}^{(k)}$ and 
coordinates $\lambda_{i,j}^{(k)}$, $1\le i < j \le 3$ and 
$1\le k \le r$.  Let $H_{[n]}$ denote the hyperplane defined by
$\sum\lambda_{i,j}^{(k)}=0$, and let $L' \subset L_2$ denote 
the set of polychrome flats.  In this notation, for $X\in L_2$, 
we have $H_X=\{\sum\lambda_{i,j}^{(k)} = 0 
\mid\lambda_{i,j}^{(k)} \in X\}$.   

The subspace $S_\Pi = H_{[n]} \cap \bigcap_{X\in L'} H_X$ arising from
the neighborly partition $\Pi$ may be realized as the nullspace of the
$r^2+1$ by $3r$ matrix $K$ whose rows are determined by $H_{[n]}$ and the
$H_X$ above.  It is readily checked that the vectors
\[
v_1 = \sum_{k=1}^r \big(e_{1,2}^{(k)} -  e_{2,3}^{(k)}\big) \quad
\text{and} \quad
v_2 = \sum_{k=1}^r \big(e_{1,3}^{(k)} -  e_{2,3}^{(k)}\big)
\]
are in $S_\Pi$.  Furthermore, it is also easy to find $3r-2$ linearly 
independent rows of the matrix $K$.  Thus $\{v_1,v_2\}$ is a basis
for $S_\Pi$.  Finally, using this basis, we see that the form 
$\langle\ ,\ \rangle_\Pi$ is trivial on $S_\Pi$.  
Thus $S_\Pi$ is an irreducible component of $\RR^1(A)$.

Since $\RR^1(A) = \cV_1(\A)$ is the tangent cone of the characteristic
variety $V_1(\A)$ at the point $\bone$, we obtain a two-dimensional 
irreducible component, $V_\Pi$, of the latter by exponentiating.
Using coordinates $t_{i,j}^{(k)}$ for the complex torus $\T$, we have
\[
V_\Pi = \biggl\lbrace
\bt \in \T
\ \biggl\vert\ 
t_{i,j}^{(1)}=\dots=t_{i,j}^{(r)}, 
1\le i< j \le 3,\ \text{and}\ 
\prod_{\genfrac{}{}{0pt}{}{1\le i < j \le 3}{1\le k\le r}}
t_{i,j}^{(k)} =1
\biggr\rbrace.
\]

If $r=3$, the set of rank one elements of the lattice of $\A=\A_{3,3,3}$ 
admits three neighborly partitions in addition to $\Pi=\Pi_1$ above 
which give rise to essential 2-tori.  
These partitions are
\begin{align*}
\Pi_2 & = \Big( 
H_{1,2}^{(1)},H_{1,3}^{(2)},H_{2,3}^{(1)} \ \Big| \ 
H_{1,2}^{(2)},H_{1,3}^{(1)},H_{2,3}^{(2)} \ \Big| \ 
H_{1,2}^{(3)},H_{1,3}^{(3)},H_{2,3}^{(3)} \Big), \\
\Pi_3 & = \Big( 
H_{1,2}^{(1)},H_{1,3}^{(1)},H_{2,3}^{(3)} \ \Big| \ 
H_{1,2}^{(2)},H_{1,3}^{(3)},H_{2,3}^{(1)} \ \Big| \ 
H_{1,2}^{(3)},H_{1,3}^{(2)},H_{2,3}^{(2)} \Big), \\
\Pi_4 & = \Big( 
H_{1,2}^{(1)},H_{1,3}^{(3)},H_{2,3}^{(2)} \ \Big| \ 
H_{1,2}^{(2)},H_{1,3}^{(2)},H_{2,3}^{(3)} \ \Big| \ 
H_{1,2}^{(3)},H_{1,3}^{(1)},H_{2,3}^{(1)} \Big).
\end{align*}
Let $V_{\Pi_k}$ denote the essential torus arising from $\Pi_k$, and 
write $V(3):= \bigcup_{k=1}^4 V_{\Pi_k}$.

\begin{lem} 
\label{lem:onetorus}
If $r \ge 4$, the variety $V(r):=V_\Pi$ is the only essential torus of 
the monomial arrangement $\A=\A_{r,r,3}$.
\end{lem}
\begin{proof}
It suffices to show that the partition $\Pi$ above is the only neighborly 
partition of $L_1(\A)$ which gives rise to an essential torus.  
Recall that $L_2(\A)$ consists of flats   
$X_{i,j}=\bigcap_{k=1}^r H_{i,j}^{(k)}$, $1\le i<j\le 3$, with $|X_{i,j}|=r$, 
and   
$Y_{p,q,s}=H_{1,2}^{(p)} \cap H_{2,3}^{(q)} \cap H_{1,3}^{(s)}$, where 
$s \equiv p+q \bmod r$, with $|Y_{p,q,s}|=3$.  

Let $\Gamma$ be a neighborly partition of $L_1(\A)$.  
If all three flats $X_{i,j}$ 
are monochrome, then either $\Gamma$ is trivial and gives rise to no essential 
torus, or $\Gamma=\Pi$ is the partition considered above, 
giving rise to the torus $V(r)$.  So assume that one of the flats $X_{i,j}$ is 
polychrome.  If all elements of $L_2(\A)$ are polychrome, 
it is readily checked that the subspace $S_\Gamma$ arising 
from $\Gamma$ consists only of the 
origin, $S_\Gamma=\{\bz\}$, and yields no essential torus.  
Thus we may assume further that there is a monochrome flat.

Without loss of generality, suppose that the flat $X_{1,2}$ is polychrome.  
If, say $X_{1,3}$ is monochrome, then since $\Gamma$ is neighborly by 
assumption, all other elements of $L_2(\A)$ must be polychrome.  
Checking that $S_\Gamma=\{\bz\}$, 
we see that $\Gamma$ does not yield an essential torus in this instance.

Suppose now that $Y_{r,r,r}$ is monochrome (and that $X_{1,2}$ is polychrome).  
If $Y_{r,r,r}$ is the only monochrome flat, one can check that 
$S_\Gamma=\{\bz\}$.  So suppose there is another monochrome flat.   
By the above considerations, we can assume that this flat is of multiplicity 
three, say $Y_{r-1,r-2,r-1}$.  
With these assumptions, the fact that $\Gamma$ is neighborly yields $3r-9$
polychrome flats: $X_{i,j}$, $1\le i<j \le 3$, 
$Y_{p,r,p}$, $Y_{p,r-p,r}$, $Y_{r,p,p}$, $1\le p\le r-1$, and 
$Y_{r-1,p+1,p}$, $Y_{p+1,r-1,p}$, $Y_{p,r-p-2,r-2}$, $1\le p\le r-3$.  
An exercise in linear algebra then reveals that 
$S_\Gamma=\{\bz\}$ in this case as well.
\end{proof}

\begin{rem} 
\label{rem:full} 
A comparable analysis reveals that the characteristic variety $V_1(\A)$ 
of every full monomial arrangement $\A=\A_{r,1,3}$ 
also contains a two-dimensional essential torus.  Briefly, we have
$\A=\A_{r,r,3} \cup\{H_1,H_2,H_3\}$, where 
$H_i=\ker x_i$, so $|\A| = n+3$.  The partition
$
\Pi = \Big(
H_3,H_{1,2}^{(i)} \ \Big| \  
H_2,H_{1,3}^{(j)} \ \Big| \  
H_1,H_{2,3}^{(k)} \Big)
$ 
of the rank one elements of $L(\A)$ is neighborly.  
All elements of $L_2(\A)$ of multiplicity greater than 2 are 
polychrome.  The subspace $S_\Pi \subset \C^{n+3}$ has basis
$\{v_1+r(e_3-e_1),v_2+r(e_2-e_1)\}$, is an irreducible component of 
$\RR^1(A)$, and gives rise to the essential torus 
\[
\bigl\{ \bt \in (\C^*)^{n+3} \ \bigl\vert\ 
\bigl(t_{i,j}^{(1)}\bigr)^r=\dots=
\bigl(t_{i,j}^{(r)}\bigr)^r = t_q, 1\le i<j \le 3, 
q\neq i,j,\ \text{and}\ T=1 \bigr\}
\]
of $V_1(\A)$, where $T=t_1t_2t_3\prod t_{i,j}^{(k)}$.
\end{rem}

Returning to the arrangement $\A=\A_{r,r,3}$, we now find all irreducible 
components of the variety $V_1(\A)$ of positive dimension.  
We first require some preliminary results.

\begin{lem} 
\label{lem:monosubarrs}
If $r=pq$, then the monomial arrangement $\A_{r,r,3}$ has $p^2$ 
subarrangements lattice-isomorphic to the arrangement $\A_{q,q,3}$.
\end{lem}
\begin{proof}
Let $K\zeta^a$ and $K\zeta^b$ be two cosets of 
$K=\langle\zeta^p\rangle=\Z_q$ in $\Z_r$.  Note that there are $p=r/q$ 
such cosets.  It is then readily checked that the subarrangement
\[
\left\{ H_{1,2}^{(a+kp)},H_{2,3}^{(b+kp)},H_{1,3}^{(a+b+kp)}
\ \bigl\vert\  
k=0,\dots,q-1 \right\},
\]
where the indices $a+b+kp$ are taken mod $r$, is lattice-isomorphic 
to $\A_{q,q,3}$.
\end{proof}

By the discussion in~\ref{sec:rr3}, each of the subarrangements $\A_{q,q,3}$ 
above gives rise to a two-dimensional component of $V_1(\A)$ (four such 
components if $q=3$).  Denote these 
components by $V(q:a,b) \cong V(q)$.  Note that the local component arising from 
$H_{1,2}^{(a)}\cap H_{2,3}^{(b)} \cap H_{1,3}^{(a+b)} \in L_2(\A)$ may
be expressed as $V(1:a,b)$, and that essential tori of
$V_1(\A)$ may be expressed as $V(r:1,1)$.  We now show that 
these components, and the local components $V_{X_{i,j}}$ arising from 
$X_{i,j} = \bigcap_{k=1}^r H_{i,j}^{(k)} \in L_2(\A)$, are the only 
positive-dimensional components of $V_1(\A)$.

\begin{lem} 
\label{lem:monsub}
If $\B$ is a subarrangement of $\A=\A_{r,r,3}$ which for each $q$ is 
not lattice-isomorphic to $\A_{q,q,3}$, then $\B$ does not give rise 
to a non-local component of $V_1(\A)$.
\end{lem}
\begin{proof}
A subarrangement $\B$ of $\A=\A_{r,r,3}$ may be specified by choosing 
three subsets $I,J,K \subset [r]$:  
$\B=\{H_{1,2}^{(i)},H_{1,3}^{(j)},H_{2,3}^{(k)} \mid i \in I, j\in 
J,k\in K\}$.  
It suffices to check that a subarrangement satisfying the above 
hypothesis admits no non-trivial 
neighborly partition.  Since $\B$ is, by assumption, not isomorphic 
to a monomial arrangement, we can find say $i \in I$ and $j\in J$ so 
that $H_{2,3}^{(i-j)}\not\in\B$.  Thus $H_{1,2}^{(i)}\cap 
H_{1,3}^{(j)} \in 
L_2(\B)$ is of multiplicity two.  If $\Pi$ is a neighborly partition 
of $\B$, then $H_{1,2}^{(i)}$ and $H_{1,3}^{(j)}$ must lie in the 
same block of $\Pi$ for this $i$ and this $j$.  It follows that 
$H_{1,2}^{(i)}$ and $H_{1,3}^{(j)}$ must lie in the same block of 
$\Pi$ for all $i$ and $j$.  
Consequently, $\Pi$ is trivial.
\end{proof}

We summarize the above discussion with the following.

\begin{prop}
\label{prop:topmono1}
The first central characteristic subvariety of the monomial arrangement 
$\A_{r,r,3}$ is given by
\begin{equation} 
\label{eqn:vrr3}
\check{V}_{1}(\A_{r,r,3}) = \bigcup_{1\le i<j \le 3} V_{X_{i,j}} \cup 
\bigcup_{q|r} V(q:a,b),
\end{equation}
and thus consists of $3$ tori of dimension $r-1$, and
$\sum_{q|r}c_q \frac{r^2}{q^2}$ tori of dimension two, 
where $c_3=4$ and $c_q=1$ for $q\neq 3$.
\end{prop}

We now use the above results to identify (non-local) components 
of the characteristic variety $V_1(\A)$ for the monomial arrangements 
$\A_{r,k,\ll}$, for general $\ll$.  
Though it is not technically a monomial arrangement, 
we first consider the braid arrangement.

\subsection{The Braid Arrangement}
\label{sec:braidarrangement}
Let $\A_\ll=\{H_{i,j} \mid 1 \le i<j \le \ll\}$ denote the braid
arrangement in $\C^\ll$, with group $P_\ll$, Artin's pure braid group.  
Let $n=\binom{\ll}{2}=|\A_\ll|$, and choose coordinates $t_{i,j}$ for
$\T$.  The rank two elements of $L(\A_\ll)$ are 
\[
L_2 = \bigl\lbrace 
H_{i,j} \cap H_{p,q}, \{i,j\}\cap \{p,q\} = \emptyset,\ 
\text{and}\ 
H_{i,j}\cap H_{i,k} \cap H_{j,k}, 1\le i<j<k\le \ll 
\bigr\rbrace.
\]
Thus the coarse combinatorial Alexander invariant is $\locB(\A_\ll) = 
\bigoplus_{1\le i<j<k\le \ll} B_{i,j,k}$, where $B_{i,j,k}$ is the 
local Alexander invariant associated to 
$H_{i,j}\cap H_{i,k} \cap H_{j,k} \in L_2$.  
As noted in Theorem~\ref{thm:Vsum}, this yields local components
$V^\loc_1(\A_\ll)=\bigcup V_{i,j,k}$ of the characteristic variety
$V_1(\A_\ll)$,
where
\[
V_{i,j,k} = \bigl\lbrace 
\bt \in \T 
\ \bigl\vert\  
t_{i,j}t_{i,k}t_{j,k}=1,\ 
\text{and}\ \:
t_{p,q} = 1\ \: \text{if}\ \:\lvert\{p,q\}\cap \{i,j,k\}\rvert \le 1
\bigr\rbrace.
\]

Non-local components of the variety $V_1(\A_\ll)$ may be
detected as follows.  Since the Coxeter groups D$_3$ and A$_3$ 
are isomorphic, the pure braid group $P_4$ coincides with the 
generalized pure braid group $P(2,2,3)$.  Thus, the calculation 
in~\ref{sec:rr3} above yields a non-local component 
$V\cong V(2)$ of $V_1(\A_4)$.  In the current notation, 
this essential torus is given by
\[
V = \bigl\lbrace 
\bt \in (\C^*)^6 \ \bigl\vert\  
t_{1,2}=t_{3,4},\ t_{1,3}=t_{2,4},\  
t_{2,3}=t_{1,4},\ {\textstyle \prod}_{i<j} t_{i,j}=1 
\bigr\rbrace.
\] 

In general, recall that $n=\binom{\ll}{2}$, write 
$T = \prod_{1\le i < j \le \ll} t_{i,j}$, and for each 
$4$-element subset $I$ of $[\ll]$, let
\begin{equation*}
V_I =  \bigl\lbrace  
\bt \in (\C^*)^n \mid 
t_{i,j}=t_{p,q}  \ \text{if}\  \{i,j\}\cup \{p,q\}= I, \ 
t_{p,q}=1 \ \text{if}\  \{p,q\}\not \subset I, \ T=1 
\bigr\rbrace.
\end{equation*} 
Each such subset $I$ gives rise to a subarrangement 
$\A_I \subset \A_\ll$ that is lattice-isomorphic to the 
arrangement $\A_4 \cong \A_{2,2,3}$.  Thus, each torus 
$V_I$ is a component of $V_1(\A_\ll)$.

\begin{prop} 
\label{prop:v1braid}
The first central characteristic subvariety of the braid 
arrangement $\A_\ll$ is given by
\[
\check{V}_{1}(\A_\ll) = V_1^\loc(\A_\ll) \cup 
\bigcup_{\genfrac{}{}{0pt}{}{I \subset [\ll]}{|I|=4}} V_I,
\]
and thus consists of $\binom{\ll}{3}+\binom{\ll}{4}=\binom{\ll+1}{4}$
tori of dimension two.
\end{prop}
\begin{proof}
The inclusion $V_1^\loc(\A_\ll) \cup \bigcup_{|I|=4} V_{I} \subseteq
\check{V}_1(\A_\ll)$ follows from the above discussion.  For the reverse
inclusion, note that every subarrangement of the braid arrangement 
is a direct product of braid arrangements of smaller rank, and that 
for $\ll\ge 5$, the arrangement $\A_\ll$ admits no non-trivial 
neighborly partition.
\end{proof}

\subsection{The Monomial Arrangements $\A_{r,r,\ll}$}
\label{sec:monoarrl}
We now obtain a similar description of the characteristic variety 
$V_1(\A)$ for the monomial arrangement $\A=\A_{r,r,\ll}$.  
Recall that the hyperplanes 
of $\A$ are denoted by $H_{i,j}^{(k)} = \ker(x_i - \zeta^k x_j)$, where 
$\zeta=\exp(2\pi\ii/r)$.
Each subset $I=\{i_1,i_2,i_3\}$ of $[\ll]$ gives
rise to a subarrangement $\A_I=\{H_{i_p,i_q}^{(k)} \mid 1\le k\le r, 
1\le p<q \le 3\}$ 
of $\A$ that is  lattice-isomorphic to $\A_{r,r,3}$.  Let 
$V_I$ denote the variety  specified by \eqref{eqn:vrr3} above, in 
the appropriate coordinates.

\begin{prop}
\label{prop:topmono2}
The first central characteristic subvariety of the monomial arrangement 
$\A_{r,r,\ll}$ is given by
\[
\check{V}_{1}(\A_{r,r,\ll}) = V^\loc_1(\A_{r,r,\ll}) \cup 
\bigcup_{\genfrac{}{}{0pt}{}{I \subset [\ll]}{|I|=3}} V_I.
\]
\end{prop}
\begin{proof}
Write $\A=\A_{r,r,\ll}$.  The inclusion 
$V_1^\loc(\A) \cup \bigcup_{|I|=3} V_I \subseteq \check{V}_{1}(\A)$  
follows from the above discussion.  For the reverse
inclusion, it suffices to show that $\check{V}_1(\A)$ 
has only those non-local components specified above.  
For this, first note that if $\ll=4$, then $\A$ has a non-trivial 
neighborly partition 
$\Pi=\Big( H_{1,2}^{(*)},H_{3,4}^{(*)} 
\ \Big|\  H_{1,3}^{(*)},H_{2,4}^{(*)} 
\ \Big|\  H_{2,3}^{(*)}, H_{1,4}^{(*)}  \Big)$.
However, an exercise in linear algebra reveals that the subspace 
$S_\Pi$ associated to this partition consists only of the origin, 
$S_\Pi=\{\bz\}$.  If $\ll \ge 5$, the arrangement $\A_{r,r,\ll}$ 
admits no non-trivial neighborly partition.  Thus, the monomial 
arrangements $\A_{r,r,\ll}$ of rank greater than three have no 
essential tori.

Finally, if $\B$ is a subarrangement of $\A=\A_{r,r,\ll}$ that is not 
lattice-isomorphic to a monomial arrangement $\A_{q,q,k}$, then an 
argument similar to the proof of Lemma~\ref{lem:monsub} shows that 
$\B$ does not give rise to a non-local component of $V_1(\A)$.
\end{proof}

\begin{rem}
\label{rem:fullmono}
The characteristic variety $V_1(\A_{r,1,\ll})$ of the full monomial 
arrangement may be analyzed in an analogous fashion, and the above 
calculations may be used to identify non-local components of this 
variety.  For instance, let $\B_\ll=\A_{2,1,\ll}$ be the Coxeter 
arrangement of type $\Bb$.  Associated to each subset $I=\{i_1,i_2,i_3\}$ 
of $[\ll]$ is a subarrangement $\B_I$ of $\B_{\ll}$ that is 
lattice-isomorphic to $\B_3=\A_{2,1,3}$.  Each such subarrangement 
gives rise to $12$ two-dimensional non-local components of $V_1(\B_\ll)$, one 
corresponding to the subarrangement itself (see Remark~\ref{rem:full}), 
and the remaining $11$ corresponding to subarrangements of $\B_I$ 
lattice-isomorphic to $\A_4$.  
In this way, we obtain $12\binom{\ll}{3}$ two-dimensional 
non-local components of $V_1(\B_\ll)$.  Similarly, for general $r$, 
subsets $I$ of $[\ll]$ as above yield subarrangements $\A_I$ of 
$\A_{r,1,\ll}$, and associated non-local components of 
$V_1(\A_{r,1,\ll})$.  Further analysis of the characteristic 
varieties of the full monomial arrangements is left to the 
interested reader.
\end{rem}

\bibliographystyle{amsalpha}

\end{document}